\documentclass[11pt,reqno]{amsart}

\usepackage{amssymb}
 
\theoremstyle{plain}
\newtheorem{theorem}{Theorem}
\newtheorem{proposition}[theorem]{Proposition}
\newtheorem{lemma}[theorem]{Lemma}
\newtheorem{corollary}[theorem]{Corollary}
\newtheorem{fact}[theorem]{Fact} 
\theoremstyle{definition}
\newtheorem{definition}[theorem]{Definition}
\theoremstyle{remark}
\newtheorem{remark}[theorem]{Remark}

\newcommand{\X}{\mathfrak X}
\newcommand{\Y}{\mathcal Y}
\newcommand{\Z}{\mathcal Z}
\newcommand{\W}{\mathcal W}

\newcommand{\e}{\varepsilon}    

\newcommand{\wt}{\widetilde}
\newcommand{\wh}{\widehat}

\newcommand{\lnm}{\left\Vert} 
\newcommand{\rnm}{\right\Vert} 
\newcommand{\lav}{\left\vert}
\newcommand{\rav}{\right\vert}
\newcommand{\elnm}{\left\vert\vert\vert} 
\newcommand{\ernm}{\vert\vert\right\vert}  
\newcommand{\lsp}{\left[}
\newcommand{\rsp}{\right]}
\newcommand{\spn}[1]{\text{\rm{sp}}\{ #1 \}}
\renewcommand{\leq}{\leqslant} 
\renewcommand{\geq}{\geqslant}

\newcommand{\bio}[2]{\left\{#1,~#2\right\}}
\newcommand{\BIO}[2]{#1~\times~#2}
\newcommand{\BIOs}[2]{#1\times#2}  
\newcommand{\ft}{fundamental total} 
\newcommand{\nc}{(2+\e)}
\newcommand{\lc}{3}  
\newcommand{\oc}[2]{#2\ominus#1} 
\newcommand{\embed}{\hookrightarrow}
\newcommand{\pperp}{\top}  
\newcommand{\sign}{\textrm{sign~}}
\newcommand{\deq}{:=}
\newcommand{\card}[1]{\left\vert #1 \right\vert}
\newcommand{\moc}{\delta} 
\newcommand{\mocw}{\moc_\W}
\newcommand{\pdz}{^*\hsp{-3}\Z_2}  

\newcommand{\hsp}[1]{\hskip #1 pt}
\newcommand{\vsp}[1]{\vskip #1 pt}
  
\begin{document}
\baselineskip 16 pt


\title[Geometry of Banach spaces and biorthogonal systems]
{Geometry of Banach spaces \\ and biorthogonal systems}
\author[Dilworth,  Girardi, Johnson]
{S.~J. Dilworth,  Maria Girardi \text{\tiny{and}}  W.~B. Johnson}  
\address{Department of Mathematics, 
University of South Carolina, \vsp{0}   
Columbia, SC 29208, U.S.A.} 
\email{dilworth@math.sc.edu and girardi@math.sc.edu}
\address{Department of Mathematics, 
Texas A\&M University,  \vsp{0}
College Station, TX 77843, U.S.A.} 
\email{johnson@math.tamu.edu}  
\thanks{Girardi is supported 
in part by NSF grant DMS-9622841 
and is a participant in the 
NSF Workshop in Linear Analysis and Probability, 
Texas A\&M  University 
(supported in part by NSF grant DMS-9311902).  
Johnson is supported in part by NSF grant 
DMS-9623260, DMS-9900185, and by Texas Advanced Research Program under
Grant No. 010366-163.} 
\subjclass{46B99, 46B25, 46B20}

\begin{abstract} 
A separable Banach space $\X$ 
contains $\ell_1$ isomorphically 
if and only if $\X$  has a 
bounded 
\ft\  $wc_0^*$-stable  biorthogonal system.
The dual of a separable Banach space~$\X$ 
fails the Schur property 
if and only if  $\X$ has a 
bounded 
\ft\  $wc_0^*$-biorthogonal system.  
\end{abstract}

\maketitle


\section{INTRODUCTION}
\label{S:intro}

Generally it is  easier to deal with  Banach spaces
that  have some sort of basis structure,  the most useful
and commonly used structures being  Schauder bases and   
finite-dimensional Schauder decompositions (FDD).  Much
research in Banach space theory has gone into proving that if a
Banach space which has  a Schauder basis or FDD possesses a certain 
property, then the space has a basis or FDD which reflects the
property.  While such theorems often give information (for example, by
passing to suitable subspaces) about general spaces which do not have
a basis or an FDD, they cannot give a classification of all separable
spaces which have a certain property  in terms of bases for the
entire space unless the property itself implies the existence of a
basis or FDD in a space which has the property.  For that reason it
is interesting to consider weaker structures than FDD's and Schauder
bases which   exist in every separable Banach space and try to prove
that a separable Banach space has a certain property if and only if
there is structure in the space which reflects the property.

One useful   basis-like structure that has been considered  for a
long time is that of fundamental total biorthogonal system.   
Markushevich~\cite{M} showed in 1943 that each separable 
Banach space contains a fundamental total biorthogonal system. 
The main theorems of this paper 
characterize   certain geometric properties 
of a Banach space by which types of 
bounded \ft\ biorthogoal systems exist in the   space.  
Theorem~\ref{t:sch}  
shows that the dual of a separable Banach space~$\X$ 
fails the Schur property 
if and only if  $\X$ contains a 
bounded  \ft\  $wc_0^*$-biorthogonal system. 
Recall that the  dual of a  Banach space~$\X$  fails 
the Schur property if and only if 
$\X$ fails the Dunford-Pettis property 
or $\ell_1$  embeds into $\X$. 
Theorem~\ref{t:elle}  
shows that  $\ell_1$  embeds in a separable Banach space $\X$
if and only if $\X$ contains a 
bounded \ft\  $wc_0^*$-stable  biorthogonal system.

Thirty-two years after  Markushevich's result~\cite[1943]{M}, 
Ovsepian and  Pe\l czy\'nski  showed~\cite{OP} that 
for each positive~$\e$,  each  separable 
Banach space contains a 
$[(1+\sqrt2)^2  + \e]$-bounded 
fundamental total biorthogonal system;
the following year  Pe\l czy\'nski~\cite{P} improved the 
bound to $(1+\e)$.   
The proofs of Theorems~\ref{t:sch} and~\ref{t:elle} 
use a combination of the methods in~\cite{OP} and ~\cite{P}. 
Theorem~\ref{t:celle} shows that 
if $\X$ is a separable Banach space containing~$\ell_1$, 
then   there is a $[1+\sqrt2  + \e]$-bounded 
fundamental  biorthogonal system
$\bio{x_n}{x_n^*}$ in $\BIO{\X}{\X^*}$ 
with    
the~ $x_n^*$'s  arbitrarily close to   an isomorphic copy of 
$\ell_2$  sitting in  $\X^*$. 
Section~\ref{S:lowerbnd}  shows that,  
in the statement of Theorem~\ref{t:celle}, 
the $(1+\sqrt{2}+\e)$ can {\bf not} 
be replaced with $(1.02 +\e)$.  
To the best of our knowledge, this is the 
first result in the literature  which provides 
the existence of a bounded fundamental biorthogonal system in
all spaces which have a certain property, and yet the bound  for the
systems cannot be arbitrarily close to one. 
  

\section{NOTATION  and TERMINOLOGY}

Throughout this paper, $\X$, $\Y$, and $\Z$  
denote arbitrary 
(infinite-dimensional real) Banach spaces.   
If $\X$ is a Banach space, 
then   
$\X^*$  is its dual space,
$B(\X)$ is its  (closed) unit ball,    
 $S(\X)$ is its unit sphere, 
$\delta \colon \X \to \X^{**}$ is the natural point-evaluation 
isometric embedding, and $\wh x = \delta (x)$.   
If $Y$ is a subset of $\X$, then 
$\spn  Y $ is the linear span of $Y$ while 
$\lsp Y \rsp$  is 
the closed linear span of~$Y$.     
Often used are 
the  unit vector basis $\{ \delta_n \}$ of $\ell_1$,   
the Kronecker delta  $\delta_{nm}$, 
and the space $C(K)$ of  
continuous functions  on a compact Hausdorff 
space~$K$.  

If $a >0$, then  $T\in\mathcal L(\X,\Y)$  
is an  \textit{$ab$--isomorphic embedding} provided   
\[
a^{-1} \lnm x \rnm ~\leq~ \lnm Tx \rnm ~\leq~ b \lnm x \rnm  
\] 
for each $x\in\X$; in this case,  
$T_o \in \mathcal L(\X, T\X)$ denotes the bijective 
operator that agrees with $T$ on $\X$.  
A surjective $\tau$--isomorphic embedding 
$T\in\mathcal L(\X,\Y)$ is a \textit{$\tau$--isomorphism};
in this case, $\X$ and $\Y$ are \textit{$\tau$--isomorphic}. 

Recall that 
for a subset $X$ of $\X$ 
and a subset $Z$  of $\X^*$  
\begin{enumerate}           
\item 
$X$  is  \textit{fundamental} if $\lsp X \rsp = \X$, 
or,  equivalently, the annihlator $X^\perp $ of $X$ in $\X^*$ is $\{ 0 \}$, 
\item 
$Z$  is  \textit{total} if the  weak$^*$-closure 
of $\spn{Z}$ is $\X^*$, 
or, equivalently,  the preannihilator $Z^\pperp$ of $Z$ 
in $\X$ is $\{ 0 \}$,  
\item
for  a fixed $\tau \geq 1$,  
$Z$ \textit{$\tau$-norms} $X$ 
(or $X$ is $\tau$-normed by $Z$)  if 
\[
\lnm x \rnm~ \leq~ \tau \, \sup_{z\in Z\setminus \{0\}} 
\frac{z(x)}{\lnm z \rnm}   
\]
for each ~$x\in X$, 
\item 
$Z$   \textit{norms} $X$ if $Z$   $1$-norms  $X$. 
\end{enumerate}
If $Z$ $\tau$-norms $\X$ for a $\tau \geq 1$      
then $Z$ is total.   
Also,  
$\bio{x_n}{x_n^*}_{n=1}^{\infty}$ 
in $\BIO{X}{Z}$ is   
\begin{enumerate} 
\item  
a   \textit{biorthogonal system}    
if $x_n^*(x_m) = \delta_{nm}$,    
\item
\textit{$M$-bounded} 
if  \, 
$\{ x_n \}$ and $\{ x^*_n \}$ are bounded 
and  $ \sup_n \lnm x_n \rnm \, \lnm x_n^* \rnm 
\leq M$, 
\item 
\textit{bounded} if it is $M$-bounded for some (finite) $M$,  
\item
\textit{fundamental} if $\{ x_n  \}$ is fundamental, 
\item 
\textit{total} if   $\{x_n^*\}$ is  total. 
\end{enumerate} 
A biorthogonal system~$\bio{x_n}{x_n^*}_{n=1}^{\infty}$ 
in~$\BIO{\X}{\X^*}$ is: 
\begin{enumerate}  
\item
a  \textit{$wc_0^*$-biorthogonal system}  
if $\{x_n^* \}$ is a semi-normalized 
(i.e., bounded and  bounded
away from zero)  weakly-null  sequence, 
\item  
a \textit{$wc_0^*$-stable biorthogonal system}  
if, for each isomorphic 
embedding $T$ of $\X $ into some  $ \Y$,   
there exists a lifting $\{ y_n^* \}$ of $\{ x_n^* \}$
(i.e., $T^* y_n^* = x_n^*$ for each $n$)
such that  
$\{ y_n^* \}$ is a semi-normalized weakly-null sequence in $\Y^*$    
(or equivalently, such that   $\bio{Tx_n}{y_n^*}$ in $\BIO{\Y}{\Y^*}$
is a $wc_0^*$-biorthogonal system). 
\end{enumerate} 
Bases of type $wc_0^*$ were introduced in~\cite{FS}   
(cf.~\cite[II.7 and pg.~625--626]{S1}).

Recall that  $\Z$ is 
\textit{injective}   if 
for each pair $\X$ and $\Y$, 
each  isomorphic  embedding $T \in\mathcal L (\X, \Y)$, 
and each $S \in \mathcal L (\X, \Z)$,  
there exists  $\wt S \in\mathcal L (\Y, \Z) $  such that 
the following diagram commutes.  
\vsp{5}
\[
\begin{picture}(100,50)
\put(10,0){\makebox(0,0){$\X$}}
\put(60,5){\makebox(0,0){\vector(1,0){50}}}
\put(60,10){\makebox(0,0){$S$}}
\put(110,0){\makebox(0,0){$\Z$}}
\put(10,25){\makebox(0,0){\vector(0,1){18}}}
\put(0,25){\makebox(0,0){$T$}}
\put(10,50){\makebox(0,0){$\Y$}}
\put(60,30){\makebox(0,0){\vector(2,-1){50}}}
\put(60,45){\makebox(0,0){$\wt S$}}
\end{picture}
\]
\vsp{5}\noindent    
If $\Z$ is injective, 
then there  exists $\lambda \geq 1$  
so that $\Z$ is \textit{$\lambda$-injective}, 
i.e.   $\wt S$ can   be chosen 
so that $\Vert \wt S \Vert \leq \lambda \lnm S T_o^{-1} \rnm$.   
Recall $\Z$ is a  \textit{Grothendieck space} if 
weak$^*$ and weak sequential convergence in $\Z^*$
coincide; an   injective space 
is a  Grothendieck space (cf.~\cite[p.~188]{LT3}). 
$\Z$ has the  \textit{Schur property} if 
weak  and  strong sequential convergence in $\Z$
coincide.

All notation and terminology, not otherwise explained, 
are as in~\cite{DU} or~\cite{LT1}.


\section{THE FINE LINE BETWEEN  
  $WC_0^*$  AND  $WC_0^*$-STABLE}

The unit vectors $\bio{e^p_n}{e^q_n}$ in 
$\BIO{\ell_p}{\ell_q}$, 
where $1 \leq p < \infty$  
and $q$ is the conjugate exponent of $p$,  
form a  $1$-bounded  
\ft\ $wc_0^*$-biorthogonal system.   
For $p=1$, they are even a 
$wc_0^*$-stable biorthogonal system, 
as the  proof of (a) implies~(b) in~Theorem~\ref{t:elle}   shows.  
The next two theorems clarify the  
fine line between the existence of 
nice $wc_0^*$-biorthogonal     
and $wc_0^*$-stable biorthogonal systems. 

\begin{theorem}
\label{t:sch}
The following statements  are equivalent. 
\begin{enumerate}
\item[\rm{(a)}]   
$\X^*$  fails the Schur property. 
\item[\rm{(b)}]   
There is a  bounded  $wc_0^*$-biorthogonal system 
             in $\BIO{\X}{\X^*}$.  
\end{enumerate}
And in the case that $\X$ is separable: 
\begin{enumerate}
\item[\rm{(c)}]   
There is a  bounded \ft\ $wc_0^*$-biorthogonal system  \\  
 $\bio{x_n}{x_n^*}$ in $\BIO{\X}{\X^*}$.  
\end{enumerate}
Furthermore for each $\e>0$:   
if {\rm (b)} holds then the system can be taken to be $(1+\e)$-bounded; 
if {\rm (c)} holds then the system   
can be taken to be $[2(1+\sqrt2)^2 + \e]$-bounded 
and so that  $\lsp x_n^* \rsp $ norms $\X$.
\end{theorem}
\noindent
Recall (cf.~\cite[p.~23]{D2}) that $\X^*$ fails 
the Schur property if and only if 
$\X$ fails the Dunford-Pettis property 
or $\ell_1 \embed \X$.
 
\begin{theorem}
\label{t:elle}
The following statements  are equivalent. 
\begin{enumerate}
\item[\rm{(a)}] $\ell_1 \embed \X$. 
\item[\rm{(b)}] There is a bounded $wc_0^*$-stable 
             biorthogonal system 
             in $\BIO{\X}{\X^*}$. 
\end{enumerate}
And in the case that $\X$ is separable:
\begin{enumerate} 
\item[\rm{(c)}] There is a bounded \ft\  $wc_0^*$-stable 
             biorthogonal system
             $\bio{x_n}{x_n^*}$ in $\BIO{\X}{\X^*}$.     
\end{enumerate} 
Furthermore for each $\e>0$:   
if {\rm (b)} holds  
then the system can be taken to be $(1+\e)$-bounded; 
if {\rm (c)} 
holds then the system can be taken to be $[(1+\sqrt2)  + \e]$-bounded 
and so that  $\lsp x_n^* \rsp $  $\nc$-norms   $\X$.  
\end{theorem} 

In this section are the proofs of the 
easier implications   in the above theorems. 
The other  implications   follow from the 
results  of the next section.

\begin{proof}[Proof of {\rm (b)} implies {\rm (a)} in~Theorem~\ref{t:sch}] 
A  $wc_0^*$-biorthogonal system 
in $\BIOs{\X}{\X^*}$ is enough to force $\X^*$ 
to fail the Schur property. 
\end{proof}

\begin{proof}[Proof of {\rm (b)} implies {\rm (a)} in~Theorem~\ref{t:elle}]   
Find an (isometric) embedding $T $ of  $\X$ into a  $C(K)$-space. 
Assume that there is a 
 $wc_0^*$-biorthogonal system 
$\bio{Tx_n}{y_n^*}$ in 
$\BIOs{C\left( K \right)}{ C^*\left( K\right)}$ 
with $\{ x_n \}$ bounded, 
which would be the case if (b) held.    
If $\{ x_n \}$ had  a weakly Cauchy 
subsequence  $\{ x_{n_k} \}$, 
then   $\{ T x_{n_k} \}$  would be weakly Cauchy 
and  $\{ y_{n_k}^*  \}$ would be weakly null, 
which cannot   
be since a ~$C(K)$ space has the Dunford-Pettis 
property  (cf.~\cite[p.~20]{D2}).  
So, by Rosenthal's $\ell_1$ theorem, 
$\{ x_n \}$ admits a subsequence that is
equivalent to the unit vector basis of~$\ell_1$. 
\end{proof}

\noindent  
The above proof reveals somewhat more.    

\begin{remark}
\label{r:iti}
In the definition of $wc_0^*$-stable biorthogonal system,   
if the word 
\textit{isomorphic} is replaced with  \textit{isometric} 
then  
the statement of Theorem~\ref{t:elle}  remains true. 
\qed
\end{remark}

\begin{remark}
\label{r:elo}
If $\bio{x_n}{x_n^*}$ in $\BIO{\X}{\X^*}$  is either:
\begin{enumerate}
\item[\rm{1.}]
a bounded $wc_0^*$-biorthogonal system 
     and $\X$ has the 
     Dunford-Pettis Property 
\end{enumerate}
or
\begin{enumerate}
\item[\rm{2.}]
a bounded $wc_0^*$-stable biorthogonal system, 
\end{enumerate} 
then each subsequence of  $ \{ x_n \}$ 
contains a further subsequence that is 
equivalent to the unit vector basis of $\ell_1$.  
\qed 
\end{remark} 

That {\rm (a)} implies {\rm (b)} 
in Theorem~\ref{t:sch} (with the $(1+\e)$ bound) 
follows from Facts~\ref{f:sch}--\ref{f:schsl}. 
 
\begin{fact}
\label{f:sch} 
Let $\{ x_n \}_{n=1}^\infty$ be a    
weakly null sequence in $\X$ and  $\{ g_n \}$ 
be a bounded sequence in $\X^*$ and $\e >0$. 
Then there exists $m\in\mathbb N$ satisfying 
$$ 
\lav \langle x_m ,  g_n \rangle \rav < \e 
$$ 
for infinitely many $n\in\mathbb N$. 
\end{fact}   
\noindent 
This follows directly from the fact that, 
since $\{ x_n \}_{n=1}^\infty$ is weakly null,  
there exists a finite sequence  $\{ \lambda_m \}_{m=1}^N$ 
of positive numbers  
satisfying  
$$ 
 \max_{\pm} \lnm \sum_{m=1}^N  \pm \lambda_m x_m \rnm  
 ~<~ \frac{\e}{ \sup_j \lnm g_j \rnm}     
$$ 
and $\sum_{m=1}^N \lambda_m = 1$  (cf.~\cite[p. 48, Exercise 13]{W}).

\begin{fact}
\label{f:schl}  
Let $\{ x_n \}_{n=1}^\infty$ be a normalized  
weakly null sequence in $\X$ and $\e >0$. 
Then there 
are a subsequence $\{ x_{n_k} \}_{k=1}^\infty$ 
and functionals $\{ x_{n_k}^* \}_{k=1}^\infty$ 
biorthogonal to $\{ x_{n_k} \}_{k=1}^\infty$ 
so that \hskip 4 pt $\sup_{k\in\mathbb N} \lnm   x_{n_k}^* \rnm < 1+ \e$.
\end{fact}

\begin{proof} 
Fix a sequence $\{ \e_k \}_{k=1}^\infty$ of positive  
numbers   
satisfying $$\sum_{k=1}^\infty \e_k < \e/6 \ . $$ 

Without loss of generality (pass to a subsequence), 
$\{ x_n \}_{n=1}^\infty$ is a basic sequence with 
biorthogonal functional $\{ f_n \}_{n=1}^\infty$ 
satisfying $\lnm f_n \rnm < 3$.  
These ~$f_n$'s will be used to perturb functionals 
as needed.    

Without loss of generality (pass to a subsequence),  
there is a    system $\bio{x_n}{g_n}_{n=1}^\infty$ 
in $\BIO{\X}{\X^*}$ satisfying  $\lnm g_n \rnm < 1 + \e/2$ and 
$$
\langle x_m , g_n  \rangle ~=~ \delta_{mn} 
\text{\hskip 20 pt when \hskip 5 pt $ n \leq m $~.}
$$   
To see how to   find such a system by induction,   
consider a subsequence 
\[ 
\{ n(j,k) \}_{k=1}^\infty 
\]
in $\mathbb N$ given at the beginning of the $j^{\text{th\/}}$ step 
(for the base step, let $n(1,k) = k$).  
Let $x_{n_j} = x_{n(j,1)}$ and     
find  $\wt{g}_{n_j}$ in $ S\left( \X^* \right) $ satisfying 
$ \wt{g}_{n_j} (x_{n_j}) = 1$.  
Find a subsequence $\{ n(j+1,k) \}_{k=1}^\infty$ 
of $\{ n(j,k) \}_{k=2}^\infty$  satisfying 
$$
\lav \langle x_{n(j+1, k)} , \wt{g}_{n_j} \rangle \rav < \e_k
$$ 
for each $k\in\mathbb N$ and let 
$$
g_{n_j} = \wt{g}_{n_j} - \sum_{k=1}^\infty 
  \langle  x_{n(j+1, k)} , \wt{g}_{n_j} \rangle  
   \  f_{n(j+1, k)}  \ . 
$$ 

Without loss of generality (pass to a subsequence),  
$$ 
\lav \langle x_m ,  g_n  \rangle \rav < \e_m   
\text{\hskip 20 pt when \hskip 5 pt $m  < n$~.}
$$ 
To    accomplish this,  
iterate Fact~\ref{f:sch} to produce a sequence 
$$\{ \{ n( j,k) \}_{k=1}^\infty \}_{j=1}^\infty$$  
of sequences and a sequence $\{ k_j \}_{j=1}^\infty$  
so that $\{ n( j+1 ,k) \}_{k=1}^\infty $ is a 
subsequence of $\{ n( j,k) \}_{k=1}^\infty $ and 
$$
\lav\langle x_{n(j, k_j)} ,  g_{n(j+1 , k)} \rangle\rav
  ~<~ \e_j \ . 
$$
Then the subsequence $n_j = n(j, k_j)$ works.  

Clearly, the functionals  
$$
  x^*_n ~\deq~ g_n ~-~ 
  \sum_{\{m \in \mathbb N \colon m < n\}} 
    \langle x_m , g_n \rangle \, f_m  
$$  
are biorthogonal to $\{ x_n \}_{n=1}^\infty$ and  
are of norm at most $1+ \e$.  
\end{proof}

\begin{fact}
\label{f:schsl} 
Let $\bio{x^*_n}{x^{**}_n}_{n=1}^\infty$ 
be a biorthogonal system in $\BIO{\X^*}{\X^{**}}$ 
with $\sup_n \ \lnm x^{**}_n \rnm < 1+\e$ for some $\e>0$
and $\{x^*_n\}$ normalized  and weak-star null.   
Then there is a subsequence $\{ n_k \}_{k=1}^\infty$ 
along with a  
biorthogonal system $\bio{x_{n_k}}{x^*_{n_k}}_{k=1}^\infty$
in $\BIO{\X}{\X^{*}}$  with 
$\sup_k \ \lnm x_{n_k} \rnm < 1+\e$.
\end{fact}

\begin{proof}   
Without loss of generality (pass to a subsequence),  
there is a biorthogonal system 
$\bio{y_n}{x^*_n}_{n=1}^\infty$
in $\BIO{\X}{\X^{*}}$  with 
$$\sup_n \ \lnm y_{n} \rnm   \leq M < \infty \ . $$
For just let $\X_0$ be a separable subspace of $\X$ 
that 1-norms $\lsp x^*_n \rsp_{n=1}^\infty$ 
and take a $\sigma(\X_0^*, \X_0)$-basic 
subsequence of $\{ x_n^*\vert_{\X_0} \}$  
(\cite{JR}, cf.~\cite[V.Exercise~7]{D1}).    

For each $n\in\mathbb N$ let 
\begin{align*} 
E_n ~&\deq~ \lsp x_n^{**} \rsp \\
F_n ~&\deq~ \lsp x_1^*, \ldots , x_n^* \rsp \ . 
\end{align*}  
Use the Principle of Local Reflexivity to 
find a sequence $\{ z_n \}_{n=1}^\infty$ in  $\X$ 
satisfying 
$$
\langle z_n , x_k^* \rangle 
~=~ 
\langle   x_k^* , x^{**}_n \rangle
~=~
\delta_{nk} 
\text{\hskip 20 pt when \hskip 5 pt $k \leq n$}
$$ 
and 
$$ 
\sup_{n\in \mathbb N} \lnm z_n \rnm < 1 + \e - \e_0 
$$ 
for some $\e_0 > 0$.    
Fix a sequence $\{ \e_j \}_{j=2}^\infty$  of 
positive  numbers satisfying 
$$
\sum_{j=2}^\infty \e_j <  
\frac{  \e_0}{M} \ . 
$$  
Without loss of generality (pass to a subsequence), 
$$
\lav \langle z_n , x^*_j \rangle \rav < \e_j 
\text{\hskip 20 pt when \hskip 5 pt  $n < j$~.} 
$$  

Clearly  the vectors  
$$
x_n \deq z_n ~-~ 
\sum_{j=n+1}^\infty \langle z_n , x_j^* \rangle \ y_j  
$$ 
are biorthogonal to $\{ x_n^* \}_{n=1}^\infty$ 
and are of norm at most  $1+\e$. 
\end{proof}

The next lemma provides a means by which  
to determine whether a $wc_0^*$-bi\-or\-thog\-o\-nal system is 
$wc_0^*$-stable.

\begin{lemma}
\label{l:fac} 
Let $\bio{x_n}{ x_n^* }$ 
be a biorthogonal system  
such that $\{ x_n^* \}$ is a semi-normalized 
       weak$^*$-null sequence in $\X^*$. 
Then $\bio{ x_n}{  x_n^* }$ 
is a $wc_0^*$-stable biorthogonal system 
if and only if 
   the operator  $S \colon \X \to c_0$
given by $$S(x) = \left( x_n^* \left( x \right)\right)$$ 
factors through an injective space.       
\end{lemma}
 
\begin{proof} 
Let $\bio{x_n}{ x_n^* }$ 
be a biorthogonal system  
such that $\{ x_n^* \}$ is a semi-normalized 
       weak$^*$-null sequence in $\X^*$. 

First, assume that the above  operator $S$  
factors through an injective space   and   
let $T \colon \X \to \Y$ be an isomorphic embedding. 
Consider the  diagram  
\vsp{5}
\[
\begin{picture}(200,50)
\put(0,0){\makebox(0,0){$\Y$}}
\put(100,0){\makebox(0,0){$\Z$}}
\put(150,5){\makebox(0,0){\vector(1,0){50}}}
\put(150,10){\makebox(0,0){$L$}}
\put(200,0){\makebox(0,0){$c_0$}}
\put(100,25){\makebox(0,0){\vector(0,-1){25}}}
\put(110,25){\makebox(0,0){$R$}}
\put(100,50){\makebox(0,0){$\X$}}
\put(150,45){\makebox(0,0){$S$}}
\put(150,30){\makebox(0,0){\vector(2,-1){45}}}
\put(50,30){\makebox(0,0){\vector(-2,-1){45}}}
\put(50,45){\makebox(0,0){$T$}}
\end{picture}
\]
\vsp{5}\noindent     
where $\Z$ is an injective space and $S = LR$.  
Since $\Z$ is injective, 
there exists ~$R_1  \in \mathcal L (\Y, \Z)$   
such that the following diagram (totally) commutes.
\vsp{5}
\[
\begin{picture}(200,50)
\put(0,0){\makebox(0,0){$\Y$}}
\put(100,0){\makebox(0,0){$\Z$}}
\put(150,5){\makebox(0,0){\vector(1,0){50}}}
\put(150,10){\makebox(0,0){$L$}}
\put(200,0){\makebox(0,0){$c_0$}}
\put(100,25){\makebox(0,0){\vector(0,-1){25}}}
\put(110,25){\makebox(0,0){$R$}}
\put(100,50){\makebox(0,0){$\X$}}
\put(150,45){\makebox(0,0){$S$}}
\put(150,30){\makebox(0,0){\vector(2,-1){45}}}
\put(50,30){\makebox(0,0){\vector(-2,-1){45}}}
\put(50,45){\makebox(0,0){$T$}}
\put(50,5){\makebox(0,0){\vector(1,0){50}}}
\put(50,10){\makebox(0,0){$R_1$}}
\end{picture}
\]
\vsp{5}\noindent 
Note that the operator $S \equiv LR$ is given by  
$$
 (LR)(x) = \left(  x^*_n \left(x\right)\right) 
\text{\quad and so \quad}
 x^*_n = R^* \, L^*  (\delta_n) \ ; 
$$
similarly, the operator $LR_1$ has the form  
$$
 (LR_1)(y) = \left( y^*_n \left(y \right)\right)  
\text{\quad where \quad}
y^*_n = R_1^*\,L^*  (\delta_n) \ . 
$$   
It is easy to check that $\bio{x_n}{x_n^*}$  
is indeed a $wc_0^*$-stable biorthogonal system:
the   commutativity  of the diagram gives 
that   $ T^* y^*_n =  x^*_n$,      
the weak-nullness 
of $\{ y^*_n \}$ follows from 
the fact that $\Z$ is a Grothendieck space 
($\{ L^* \delta_n \}$ is  weak$^*$-null and thus weakly-null), 
and $\lnm  y^*_n \rnm 
\geq \lnm T \rnm^{-1} \lnm x_n^*\rnm$.

Next assume that 
$\bio{x_n}{x_n^*}$  
is   a $wc_0^*$-stable biorthogonal system.  
Find an embedding $R$ 
from $\X$  into the injective space ~ $\ell_\infty(\Gamma)$ 
for   some index set~ $\Gamma$.  
By the stability of the system, 
there exists a weakly-null sequence 
$\{ y_n^* \}$ in ~$\ell^*_\infty(\Gamma)$ such 
that $R^*y^*_n = x^*_n$.    
Define $L \colon \ell_\infty(\Gamma) \to c_0$ by
$L(f)   =  \left( y^*_n \left( f \right)\right)$, 
for then  ~ $S = LR$.  
\end{proof}

\noindent
The  commutative diagram in the 
next proof   was inspired  by the 
Hagler--Johnson proof~\cite{HJ} 
of the Josefson and 
Nissenzweig Theorem 
(cf.~\cite[Chapter~XII]{D1}).  

\begin{proof}[Proof of {\rm (a)} implies {\rm (b)}  
in Theorem~\ref{t:elle}, along with the $(1 + \e)$ bound]
\hfill 
Consider the     
following commutative diagram   
\vsp{5}
\[
\begin{picture}(200,50)
\put(0,0){\makebox(0,0){$\X$}}
\put(100,0){\makebox(0,0){$L_\infty$}}
\put(150,5){\makebox(0,0){\vector(1,0){50}}}
\put(150,10){\makebox(0,0){$L$}}
\put(200,0){\makebox(0,0){$c_0$}}
\put(100,25){\makebox(0,0){\vector(0,-1){25}}}
\put(110,25){\makebox(0,0){$R_0$}}
\put(100,50){\makebox(0,0){$\ell_1$}}
\put(150,45){\makebox(0,0){$i$}}
\put(150,30){\makebox(0,0){\vector(2,-1){45}}}
\put(50,30){\makebox(0,0){\vector(-2,-1){45}}}
\put(50,45){\makebox(0,0){$j$}}
\end{picture}
\]
\vsp{5}\noindent 
where  $j$ is an isomorphic embedding,    
$i$ is  the  formal injection, 
and 
$$
R_0 (\delta_n) = r_n
\text{\qquad and\qquad}
L(f) = \left( \int f r_n \, d\mu  \right)_n
$$    
for the Rademacher functions $\{  r_n \}$. 
Since $L_\infty$ is $1$-injective, 
there exists an  operator
$R_2$ such that  the following  diagram   commutes           
\vsp{5}
\[
\begin{picture}(200,50)
\put(0,0){\makebox(0,0){$\X$}}
\put(100,0){\makebox(0,0){$L_\infty$}}
\put(150,5){\makebox(0,0){\vector(1,0){50}}}
\put(150,10){\makebox(0,0){$L$}}
\put(200,0){\makebox(0,0){$c_0$}}
\put(100,25){\makebox(0,0){\vector(0,-1){25}}}
\put(110,25){\makebox(0,0){$R_0$}}
\put(100,50){\makebox(0,0){$\ell_1$}}
\put(150,45){\makebox(0,0){$i$}}
\put(150,30){\makebox(0,0){\vector(2,-1){45}}}
\put(50,30){\makebox(0,0){\vector(-2,-1){45}}}
\put(50,45){\makebox(0,0){$j$}}
\put(50,5){\makebox(0,0){\vector(1,0){50}}}
\put(50,10){\makebox(0,0){$R_2$}}
\end{picture}
\]
\vsp{5}\noindent 
and  $\lnm R_2 \rnm \leq \lnm R_0 j_o^{-1} \rnm$.

The operator $S \deq L\,R_2$ takes the form 
$$
S(x) ~=~ \left( x^*_n \left( x \right)\right) 
\text{\quad where \quad} 
  x_n^* =  R_2^* \, L^*  \, (\delta_n) \ . 
$$
It is easy to check that 
 $\bio{j \delta_n}{x_n^*}$ is a $wc_0^*$-stable biorthogonal system.    
Biorthogonality follows from the commutativity  of the diagram.  
Since $\{ \delta_n \}$ is weak$^*$-null, so is $\{ x_n^* \}$.    
$L_\infty$ is an injective space through which $S$ factors. 
Furthermore, since $R_0$ and $L$   both have  norm one and 
$1 = x^*_n (j\delta_n)$,    
$$
\lnm j \rnm^{-1} 
~\leq~ 
\lnm x_n^* \rnm 
~=~  
\lnm R^*_2 L^* \delta_n  \rnm
~\leq~ 
\lnm R_2 \rnm \lnm L \rnm 
~\leq~ 
\lnm R_0 j_o^{-1} \rnm 
~\leq~ 
\lnm j^{-1}_o \rnm   
$$ 
and so $\{ x^*_n \}$ is semi-normalized. 
If $\ell_1$ embeds into $\X$, 
then  it $(1+\e)$-embeds  into~$\X$, 
thus one can arrange that 
$ 
\sup_n \lnm x_n \rnm \lnm x_n^* \rnm 
\leq 
\lnm j \rnm \lnm j^{-1}_o \rnm  
\leq 1 + \e 
$~.
\end{proof}

It is not difficult to see that, 
if $\X$ is any Banach space and $\e >0$, 
then there is a $(2+\e)$-bounded 
biorthogonal system $\bio{x_n}{x^*_n}$ 
in $\BIO{\X}{\X^*}$ with 
$\{ x_n^* \}$ weak$^*$-null.  
The first step  towards this is the lemma below.

\begin{lemma}
\label{l:fcd} 
If $\X_0$ is a finite co-dimensional subspace of $\X$ 
and $\e>0$, 
then there is a weak$^*$-closed finite co-dimensional 
subspace $\Y$ of $\X^*$
such that $\Y$ is  $(2 + \e)$-normed by $\X_0$. 
\end{lemma} 

\noindent 
To see how to use  Lemma~\ref{l:fcd} to produce 
the desired biorthogonal system $\bio{x_n}{x^*_n}_{n=1}^\infty$, 
start with a normalized weak$^*$-null 
sequence $\{ y_n^* \}_{n=1}^\infty$ in $\X^*$ 
(guaranteed to exist by 
the Josefson-Nissenzweig Theorem) 
and fix a sequence $\{ \e_n \}_{n=1}^\infty$ 
of positive numbers tending to zero.  
Assume that  $$\bio{x_j}{x^*_j}_{j< n }$$ have 
been found.  Let 
$$
\X_n =\lsp x_j^*\rsp_{j < n}^\pperp
\text{\quad and \quad}
 \Z_n =\lsp x_j\rsp_{j < n}^\perp \ . 
$$
By Lemma~\ref{l:fcd}, there is  a 
weak$^*$-closed  finite co-dimensional subspace $\Y_n$ of ~$\X^*$ 
that is $(2+ \e/2 )$-normed by $\X_n$. 
Since $\Y_n \cap \Z_n$ is finite co-dimensional 
and weak$^*$-closed and  
$\{ y_n^* \}_{n=1}^\infty$ is weak$^*$-null, 
there exists $x^*_n \in S\left(\Y_n \cap \Z_n\right)$ 
with $\lnm x^*_n - y^*_{k_n} \rnm < \e_n$ 
for some large~$k_n$. 
Next find $\wt{x}_n \in S\left( \X_n \right)$ with 
$1 \leq (2 + \e) x^*_n \left( \wt{x}_n \right)$ 
and let $x_n \deq  \wt{x}_n / x^*_n \left( \wt{x}_n \right)$.

\begin{proof}[Proof of Lemma~\ref{l:fcd}]  
Let $\Y$ be the annihilator of any 
finite dimensional subspace of~$\X$ that   
$(1+\e)$-norms the annihilator of $\X_0$.
For then if $f \in S( \Y)$ then 
\begin{align*}
\sup_{x_0\in S(\X_0)} \, \lav f(x_0) \rav ~
&=~ 
\inf_{y^* \in \X_0^\perp} \lnm f - y^* \rnm \\
&\geq~ 
\inf_{y^* \in \X_0^\perp}  ~\max~  
\big[ \lnm f \rnm - \lnm y^* \rnm \ , \ 
    \sup_{x \in S(\Y^\pperp)}  \lav \left( f- y^* \right) \left(x \right)\rav \big] \\
&\geq~ 
\inf_{y^* \in \X_0^\perp}  ~\max~  
\big[ 1 - \lnm y^* \rnm \ , 
\ 
     \tfrac{1}{1+\e} \lnm y^* \rnm  \big] \\ 
&\geq~ 
\inf_{0 \leq t < \infty}  ~\max~ 
\big[ 1 - t \ , \  \tfrac{t}{1+\e}\big] \\
&=~  \left( 2+\e \right)^{-1}   \  .     
\end{align*}
\end{proof} 

Lemma~\ref{l:fcd}  is  nearly  best possible  
since, for each $\e>0$,  the 
one co-dimensional subspace $\X_0$ of mean zero functions 
in $L_1$ does not $(2-\e)$-norm any finite 
co-dimensional subspace   of $L_\infty$.  
Indeed, any finite co-dimensional subspace of $L_\infty$ 
contains a norm one functional $y^*$ 
that is bounded below by $-\e$  
(just perturb 
a disjointly supported sequence of nonnegative norm one functions in
$L_\infty$   that are close to $\Y$)     
and so  $y^*(x)   \leq \frac12\left(1+\e\right)$ 
for each $x\in S\left(\X_0\right)$.  
However, any one co-dimensional subspace $\X_0$ of 
a Banach space~$\X$ does $2$-norm a  
one co-dimensional subspace 
$\Y$, 
namely $\Y \deq \text{ker\,}P$ where  
$P\colon \X^* \to \X_0^\perp$ is   a norm one projection. 
Indeed,   if $f\in S\left(\Y\right)$ then 
\begin{align*} 
\lnm f - y^* \rnm 
~&\geq~ 
\frac12~ 
\left[~ \lnm f - y^* \rnm ~+~ \lnm P\left(f - y^*\right) \rnm~\right] \\
~&=~
\frac12~ 
\left[~ \lnm f - y^* \rnm ~+~ \lnm  y^*  \rnm ~\right] 
~\geq~  
\frac{1}{2}     
\end{align*}
for each $y^* \in \X_0^\perp$.


\section{CONSTRUCTING FUNDAMENTAL TOTAL $wc_0^*$-BIORTHOGONAL SYSTEMS}

The constructions of  fundamental total  
biorthogonal systems in the 
proofs of~(a) implies~(c) in  Theorems~\ref{t:sch} and ~\ref{t:elle}  
use  the Haar matrices,  
which are summarized below.

\begin{remark}
\label{r:hwm}   
Fix $m \geq 0$ and consider the 
$2^m$-dimensional Hilbert space $\ell_2^{2^m}$, 
along with its unit vector basis $\{ e^2_j \}_{j=1}^{2^m}$. 

The   Haar basis $\{ h^m_j \}_{j=1}^{2^m}$ 
of $\ell_2^{2^m}$ can be described as follows.  For $0 \leq n \leq m$ 
and $1 \leq k \leq 2^n$  let 
\[
 I^n_k ~=~ \left\{ j \in \mathbb N \hsp{3}\colon \hsp{3}
2^{m-n} \, (k-1) \hsp{3}<\hsp{3} 
j \hsp{3}\leq\hsp{3} 
2^{m-n} \,  k  \right\}\ . 
\]  
Thus 
\begin{gather*}
I^0_1 = \left\{ 1, 2 \hsp{2},\hsp{2}  \ldots \hsp{2},\hsp{2} 2^m \right\} \\
I^1_1  = \left\{ 1, 2 \hsp{2},\hsp{2}\ldots \hsp{2},\hsp{2} 
          2^{m-1} \right\}
\textrm{\hsp{10}and\hsp{10}}
I^1_1  = \left\{ 1 + 2^{m-1}  \hsp{2},\hsp{2} \ldots \hsp{2},\hsp{2}  
            2^m \right\} \ . 
\end{gather*}
 In general, 
the collection $\{ I^n_k \}_{k=1}^{2^n}$ 
of sets along the  $n^{\text{th\/}}$-level (disjointly) partitions 
$\{ 1, 2, \ldots , 2^m \}$ into $2^n$ sets, each containing 
$2^{m-n}$ consecutive integers, 
and $I^n_k$ is the disjoint union $I^n_k = I^{n+1}_{2k-1}  \cup I^{n+1}_{2k} $.  
Now let 
\[ 
h^m_1 ~=~ 2^{\frac{-m}{2}} ~ \sum_{j\in I^0_1} e^2_j 
\] 
and, for $0 \leq n < m$ and $1 \leq k \leq 2^n$, 
let $h^m_{2^n + k}$ be supported on $I^n_k$ as 
\[
h^m_{2^n + k} ~=~ 2^{\frac{n-m}{2}} ~ 
\left[ \sum_{j\in I^{n+1}_{2k-1}} e^2_j  
~-~ \sum_{j\in I^{n+1}_{2k}} e^2_j  \right] \ . 
\] 
Note that   $\{ h^m_j \}_{j=1}^{2^m}$ 
forms an orthonormal basis for~$\ell_2^{2^m}$. 
    
Let $H_m = \left( a^m_{ij} \right)$ be the $2^m \times 2^m$ 
Haar matrix  that transforms the unit vector basis  of $\ell_2^{2^m}$  
onto the Haar basis; thus, the $j^{\text{th\,}}$ column vector of 
$H_m$ is just  $h^m_{j}$ 
and so $H_m$ is a unitary matrix.   
For example, 
$$
H_2 ~=~  \bmatrix 
2^{-1}  &  +2^{-1}  & +2^{-1/2}  & 0 \\
2^{-1}  &  +2^{-1}  & -2^{-1/2}  & 0 \\
2^{-1}  &  -2^{-1} &  0          & +2^{-1/2}  \\
2^{-1}  &  -2^{-1} &  0          & -2^{-1/2}  
\endbmatrix ~ .
$$

Let $\bio{z_j}{z_j^*}_{j=1}^{2^m}$  be a 
biorthogonal sequence in $\BIO{S(\X)}{\X^*}$. 
Consider $\bio{x_i}{x_i^*}_{i=1}^{2^m}$ where
\begin{align*}
H_m ~ \bmatrix  z_1 \\ \vdots \\  z_{2^m} \endbmatrix 
~=~ \bmatrix  x_1 \\ \vdots \\  x_{2^m} \endbmatrix  
\hskip 40 pt &\text{and}\hskip 40 pt    
H_m ~ \bmatrix  z_1^* \\ \vdots \\  z_{2^m}^* \endbmatrix 
~=~ \bmatrix   x_1^* \\ \vdots \\   x_{2^m}^* \endbmatrix \ , \\
\intertext{thus}
x_i \deq \sum_{j=1}^{2^m} a^m_{ij}  z_{j} 
\hskip 40 pt &\text{and}\hskip 40 pt       
x_i^*  \deq \sum_{j=1}^{2^m} a^m_{ij}  z_{j}^*   \ . 
\end{align*}   
Since $H_m $ is a unitary  matrix   
\begin{enumerate} 
\item[\rm{(H1)}]   
$x_i^* (x_j) = \delta_{ij} $
\item[\rm{(H2)}]  $\lsp x_i \rsp_{i=1}^{2^m} 
        = \lsp   z_j \rsp_{j=1}^{2^m} $
\item[\rm{(H3)}]  $\lsp x_i^* \rsp_{i=1}^{2^m} 
       = \lsp   z_j^* \rsp_{j=1}^{2^m} $.
\end{enumerate}
Note that, for each $1\leq i\leq 2^m$, 
\begin{enumerate}
\item[\rm{(H4)}]  $ a_{i1}^m = 2^{{-m}/{2}} $ 
\end{enumerate} 
and 
the $\ell_1$-norm of the $i^{\text{th\/}}$ row of $H_m$ is bounded  
\begin{enumerate} 
\item[\rm{(H5)}]  $\sum _{j=1}^{2^m} \lav a^m_{ij}\rav 
            ~=~ 1 + \sqrt 2 - 2^{\frac{1-m}{2}} 
              ~  \ ^{\ _{m \to \infty}} \hskip - 8 pt \nearrow  
              ~ 1 + \sqrt 2 $   
\end{enumerate} 
and so  
\begin{enumerate} 
\item[\rm{(H6)}] $\lnm   x_i \rnm \ \ ~ \leq ~
     2^{{-m}/{2}} \lnm  z_{1}  \rnm \   ~+~ 
      \left( 1 + \sqrt 2  \right)
           \max_{ 1 < j \leq 2^m }   \lnm   z_j  \rnm$  
\item[\rm{(H7)}]  $\lnm  x_i^* \rnm ~ \leq ~ 
       2^{{-m}/{2}} \lnm  z_{1}^* \rnm   ~+~ 
      \left( 1 + \sqrt 2  \right)
           \max_{ 1 < j \leq 2^m }   \lnm   z_j^* \rnm$  
\item[\rm{(H8)}] for each $x^{**} \in \X^{**}$  \newline 
  $\lav x^{**} (x_i^*) \rav    ~\leq~  $     
    $\lnm  x^{**} \rnm 2^{{-m}/{2}}  \lnm    z_{1}^* \rnm  ~+~
           \left( 1 + \sqrt 2  \right)
          \max_{ 1 < j \leq  2^m } 
         \lav  x^{**} \left(   z_j^*  \right) \rav $ . 
\end{enumerate}  
\end{remark}

\begin{definition}
\label{d:blk} A sequence $\{ J_k \}_{k=1}^\infty$ of 
subsets of $\mathbb N$ is a  \textit{blocking} 
of $\mathbb N$ if $\mathbb N$ is the 
disjoint union $\cup_{k=1}^\infty J_k$ and 
$$
   \max J_k < \min J_{k+1}
$$ 
for each $k\in \mathbb N$.  Given a blocking  $\{ J_k \}_{k=1}^\infty$ 
of $\mathbb N$, let $J_0 = \{ 0 \}$ and 
\begin{gather*} 
J_k^p   ~\deq~  \bigcup\limits_{0\leq j<k} J_j \\
J_k^o  ~\deq~ J_k\setminus\{ \text{the first element in } J_k\}  \\
J_k^{po}  ~\deq~  \bigcup\limits_{0\leq j<k}  J_k^o \\ 
\mathbb N^o ~\deq~ \bigcup\limits_{k=1}^\infty  J_k^o
\end{gather*}  
for each $k \in \mathbb N$. 
\end{definition}

\hsp{0} From the next theorem it easily follows, when ~$\X$ is separable,  
that~(a) implies~(c) in~Theorem~\ref{t:sch}.

\begin{theorem}
\label{t:csch}
Let $\X^*$ fail the Schur property. 
Fix $\e > 0$ 
along with 
$\bio{ a_n }{ b_n^*}$ in
$\BIO{\X}{\X^*}$.    
Then there exists a  
 $[ 2 (1+\sqrt{2})^2 + \e]$-bounded
$wc_0^*$-biorthogonal system 
$\bio{x_n}{x_n^*}$ in $\BIO{\X}{\X^*}$ 
such that 
$\lsp a_n \rsp \subset \lsp x_n \rsp  $ 
and $\lsp b_n^*  \rsp \subset \lsp x_n^* \rsp  $.
\end{theorem}

\begin{proof}  
Without loss of generality, 
$[ a_n ]_{n\in\mathbb N}$ and  $[ b^*_n ]_{n\in\mathbb N}$   
are each infinite dimensional.   
Fix  a sequence $\{ \delta_k \}_{k=1}^\infty$ 
of positive numbers decreasing to zero.  
Since $\X^*$ fails the Schur property, 
there is a weakly-null   
sequence ~$\{ w_i^*\}_{i=1}^\infty$ in $S(\X^*)$.

It suffices to find a system 
$\bio{x_n}{x_n^*}_{n=1}^\infty$ in $\BIO{\X}{\X^*}$  
along with 
a blocking $\{ J_k \}_{k=1}^\infty$ of~$\mathbb N$,
a sequence~ $\{ \beta_n \}_{n\in \mathbb N^o}$  from  $(0, 2+\e]$, 
and 
an increasing sequence $\{ i_n\}_{n\in \mathbb N^o}$ from $\mathbb N$,  
satisfying   
\begin{enumerate}
\item[\rm{(1)}]   
 $x_m^*(x_n) = \delta_{mn}$ 
\item[\rm{(2)}]   
 $\lnm x_n \rnm ~\leq~ \left(1+\sqrt{2} \right) + \e$
\item[\rm{(3)}]   
 $\lnm x_n^* \rnm ~\leq~  \left(2+\e\right)\,\left(1+\sqrt{2} \right) + \e $
\item[\rm{(4)}]   
 for each $x^{**}\in S\left(X^{**}\right)$, if $n\in J_k$ then 
\newline
$\lav x^{**}\left(x_n^*\right) \rav ~\leq~
\delta_k + \left( 1+\sqrt2\right) \, 
\max_{j \in J^o_k} 
\left(~ \lav x^{**}\left( \beta_j w_{i_j}^* \right) \rav 
      + \beta_j \delta_k~\right)$ 
\item[\rm{(5)}]   
 $\lsp a_n \rsp_{n=1}^{\infty} 
       \subset \lsp x_n \rsp_{n=1}^{\infty}$
\item[\rm{(6)}]   
  $\lsp b^*_n \rsp_{n=1}^{\infty} 
       \subset \lsp x_n^* \rsp_{n=1}^{\infty} $   \ . 
\end{enumerate}

The construction will inductively produce 
blocks $\bio{x_n}{x_n^*}_{n\in J_k}$.  
Let $x_0$ and $x^*_0$ be the zero vectors  
and $j_0 = 0$.  
Fix $k \geq 1$.   
Assume that $\{ J_j \}_{0\leq j<k}$ 
along with 
$\bio{x_n}{x_n^*}_{n \in J_k^p}$ 
and  $\{ i_n \}_{n \in J_k^{po}}$ 
and  $\{ \beta_n \}_{n\in J_k^{po}}$   
have been constructed to satisfy conditions~(1) through (4).   
Now to construct $J_k$ 
along with $\bio{x_n}{x_n^*}_{n \in J_k }$ 
and   $\{ i_n \}_{n  \in J_k^o }$ and $\{ \beta_n \}_{n  \in J_k^o }$. 

Let 
$$
 \mathcal P_k ~\deq~ 
    \left[ x^*_n \right]_{n\in J_k^p }^\pperp  
    \text{\hskip 20 pt and \hskip 20 pt} 
 \mathcal  Q_k ~\deq~ 
    \left[ x_n \right]_{n \in J_k^p}^\perp   
$$  
and 
$$
    n_k = \max  J^p_k \ . 
$$
The idea is to    find a biorthogonal system 
$\bio{z_n}{z_n^*}_{n\in J_k}$ in $\BIO{\mathcal  P_k}{\mathcal  Q_k}$ 
by first finding  $\{z_{1+n_k}, z^*_{1+n_k} \}$ which helps 
guarantee condition~(5) if $k$ is odd  and     
condition~(6) if $k$   even; 
however,   $\{z_{1+n_k}, z^*_{1+n_k} \}$ 
would not necessarily satisfy conditions~(2) through~(4) 
and so $J^o_k$ and 
$$
    \bio{z_n}{z_n^*}_{n \in J_k^o} \ , 
$$  
along with $\{ i_n \}_{n\in J^o_k}$  
and $\{ \beta_n \}_{n\in J^o_k}$
are constructed and then the  Haar matrix is applied to 
$\bio{z_n}{z_n^*}_{n\in J_k}$ to produce $\bio{x_n}{x_n^*}_{n\in J_k}$  
so that   
$$
  \bio{x_n}{x_n^*}_{n\in  J^p_k \, \cup J_k }
$$
with $\{ i_n \}_{n\in J^{po}_k \cup J^o_k}$  
and $\{ \beta_n \}_{n\in J^{po}_k \cup J^o_k}$ 
satisfy conditions~(1) through (4).

$\bio{z_{1+n_k}}{z^*_{1+n_k}}$ is constructed 
by a standard Gram-Schmidt biorthogonal procedure. 
If $k$ is odd, start  in~$\X$. Let 
$$ 
h_k ~=~ \min \left\{ h \colon 
     a_h \not\in \lsp x_n \rsp_{n \leq n_k} \right\} \ .  
$$ 
Set 
\begin{align*} 
z_{1+n_k} 
~&=~ 
a_{h_k} ~-~  
   \sum_{n \leq n_k} x_n^*(a_{h_k}) x_n \ ,  \\ 
\intertext{and for any $y^*_{1+n_k}$ in $\X^*$ such 
that   $y^*_{1+n_k}(z_{1+n_k}) \neq 0$,}  
z^*_{1+n_k}
~&=~ 
\frac{
y^*_{1+n_k}  ~-~ 
   \sum_{n \leq n_k}  y^*_{1+n_k}  (x_n) x_n^*}
{  y^*_{1+n_k}  (z_{1+n_k})} \ . 
\end{align*}  
If $k$ is even, start  in~$\X^*$. Let 
$$ 
h_{k} ~=~ \min \left\{ h \colon 
     b^*_h \not\in \lsp x_n^* \rsp_{n\leq n_k} \right\} \ .  
$$ 
Set 
\begin{align*} 
z_{1+n_k}^* 
~&=~ 
b^*_{h_k} ~-~  
   \sum_{n \leq n_k} b^*_{h_k} (x_n) x_n^* \ ,  \\ 
\intertext{and, for any $y_{1+n_k}$ in $\X$  such 
that   $ z_{1+n_k}^* (y_{1+n_k})  \neq 0$,}  
z_{1+n_k}
~&=~ 
\frac{
y_{1+n_k}  ~-~ 
   \sum_{n\leq n_k}  x^*_n ( y_{1+n_k} )  x_n}
{ z_{1+n_k}^* (y_{1+n_k}) } \ . 
\end{align*}
Clearly  $z_{1+n_k}^* \left( z_{1+n_k} \right) = 1$ 
and 
$$
  z_{1+n_k} ~ \in ~ \mathcal  P_k 
  \text{\hskip 25 pt and \hskip 25 pt}
    z_{1+n_k}^* ~ \in ~ \mathcal  Q_k   \ . 
$$

Find a natural number  $m_k$ larger than one  so that 
$$
2^{{-m_k}/{2}} \ 
\max\left(~ \lnm z_{1+n_k} \rnm ~,~ \lnm z_{1+n_k}^* \rnm ~\right) 
~<~ 
\min\left(~ \e ~,~ \delta_k ~\right)    
$$
and let 
$$
  J_k \deq \{ 1 + n_k , \ldots ,  2^{m_k} + n_k \} 
\text{\hskip 10 pt and  so \hskip 10 pt} 
  J_k^o \deq \{ 2 + n_k , \ldots , 2^{m_k} + n_k \} \ . 
$$  

Let
$$
\wt{\mathcal  P}_k ~\deq~  \mathcal  P_k \cap \left[ z^*_{1+n_k}\right]^\pperp
\text{\hskip 25 pt and \hskip 25 pt} 
\wt{\mathcal  Q}_k ~\deq~ \mathcal  Q_k \cap \left[ z_{1+n_k}\right]^\perp \ . 
$$  
The next step is to   find a biorthogonal system  
$\bio{z_n}{z_n^*}_{n\in J_k^o}$ 
along with $\{ i_n \}_{n\in J_k^o}$ 
and $\{ \beta_n \}_{n\in J_k^o}$ 
satisfying 
\begin{equation}
\label{tagz}
\bio{z_n}{z_n^*}  \in 
  \BIO{S\left(\wt{\mathcal  P}_k \right)}
  {\left(2+\e\right) B\left(\wt{\mathcal  Q}_k \right)}
\end{equation}
and 
\begin{equation}
\label{tagzz}
\lnm \frac{z^*_n}{\beta_{n}} - w^*_{i_n} \rnm  ~<~ \delta_k 
\end{equation} 
for each $n\in J_k^o$.
Towards this, 
fix $j\in J_k^o$   and assume that 
a biorthogonal system   
$\bio{z_n}{z_n^*}_{2 + n_k \leq n < j}$   
along with     $\{ i_n \}_{  2 + n_k \leq n < j  } $ 
and   $\{ \beta_n \}_{  2 + n_k \leq n < j  } $ 
have been constructed so that conditions~\eqref{tagz} 
and~\eqref{tagzz}   
hold for $2 + n_k \leq n < j$. 
Let 
$$
\X_{j} ~\deq~ 
\wt{\mathcal  P}_k \cap \left[ z^*_n \right]_{2+n_k \leq n < j  }^\pperp
\text{\hskip  25 pt and \hskip  25 pt} 
\Y_{j} ~\deq~ 
\wt{\mathcal  Q}_k \cap \left[ z_n \right]_{2+n_k \leq n < j  }^\perp \ . 
$$
By~Lemma~\ref{l:fcd}, there is a weak$^*$-closed 
finite co-dimensional subspace 
$\wt \Y_{j} $ of $\X^*$ such that $\wt \Y_{j}$ is 
$(2+ \e/2)$-normed by $\X_{j}$. 
Find $i_{j} > i_{j - 1}$ 
and  $ y^*_{j} \in S\left( \Y_{j} \cap \wt \Y_{j}  \right)$ 
such that 
$$  
    \lnm   y^*_{j} - w^*_{i_{j}} \rnm < \delta_k \ . 
$$  
Find $z_{j} \in S(\X_{j})$ 
such that 
$$
\frac{1}{2+\e} ~\leq~   y^*_{j} \left(z_{j} \right) 
      ~\deq~ \frac{1}{\beta_{j}} 
$$ 
and normalize 
$$
 z^*_{j}  ~\deq~  \beta_{j} ~ y^*_{j}  \ . 
$$
This completes the inductive construction of 
$\bio{z_n}{z_n^*}_{n\in J_k^o}$  
along with the sets  
$\{ i_n \}_{n\in J_k^o}$ and $\{ \beta_n  \}_{n\in J_k^o}$.

Now apply the Haar matrix  
to $\bio{z_n}{z_n^*}_{n\in J_k}$ to 
produce  $\bio{x_n}{x_n^*}_{n\in J_k}$.    
With help from the observations in 
Remark~\ref{r:hwm}, note that 
$\bio{x_n}{x_n^*}_{n\in J_k}$ 
is biorthogonal and is in $\BIO{\mathcal  P_k}{\mathcal  Q_k}$. 
Furthermore, 
for each $n$ in $J_k$,   
\begin{align*} 
\lnm x_n \rnm 
~&\leq~ 
2^{{-m_k}/{2}} \   \lnm z_{1+n_k} \rnm ~+~ 
\left( 1+ \sqrt{2} \right) \max_{j\in J^o_k} \lnm z_{j} \rnm \\
~&\leq~ 
\e ~+~ \left( 1+ \sqrt{2} \right)
\end{align*}   
and 
\begin{align*} 
\lnm x_n^* \rnm 
~&\leq~ 
2^{{-m_k}/{2}} \   \lnm z_{1+n_k}^* \rnm ~+~ 
\left( 1+ \sqrt{2} \right) \max_{j\in J^o_k} \lnm z_{j}^* \rnm \\
~&\leq~ 
\e ~+~ \left(2+\e\right)\,\left( 1+ \sqrt{2} \right)
\end{align*}  
and for each $x^{**}\in S\left(\X^{**}\right)$ 
\begin{align*} 
\lav x^{**}\left(x_n^*\right) \rav 
~&\leq~ 
2^{{-m_k}/{2}} \   \lnm z_{1+n_k}^* \rnm ~+~ 
\left( 1+ \sqrt{2} \right) 
\max_{j\in J^o_k} \lav x^{**}\left(z_{j}^*\right) \rav \\
~&\leq~ 
\delta_k ~+~  \left( 1+ \sqrt{2} \right)
\max_{j\in J^o_k} \left(~ \lav x^{**}\left(\beta_j w_{i_j}^*\right) \rav
   + \lav \beta_j \delta_k \rav~ \right) \ . 
\end{align*}

Thus 
$$
  \bio{x_n}{x_n^*}_{n\in  J^p_k \, \cup J_k }
$$
with $\{ i_n \}_{n\in J^{po}_k \cup J^o_k}$  
and $\{ \beta_n \}_{n\in J^{po}_k \cup J^o_k}$ 
satisfy conditions~(1) through (4).  
If $k$ is odd, then 
$$
\left[a_h\right]_{h\leq h_k}  \in \lsp x_n , z_{1+n_k} \rsp_{n\in J^p_k} 
        \subset \lsp x_n   \rsp_{n\in J_k^p\cup J_k} \ , 
$$
while if $k$ is even, then 
$$
\left[b^*_h\right]_{h\leq h_k} \in \lsp x_n^* , z_{1+n_k}^* \rsp_{n\in J^p_k} 
        \subset \lsp x_n^*   \rsp_{n\in J_k^p\cup J_k} \ .  
$$

Clearly the constructed system    
$\bio{x_n}{x_n^*}_{n=1}^\infty$, 
with  the  blocking $\{ J_k \}_{k=1}^\infty$ of $\mathbb N$ and 
the increasing sequence $\{ i_n\}_{n\in \mathbb N^o}$ from $\mathbb N$, 
and the sequence $\{ \beta_n\}_{\in \mathbb N^o}$ 
from $\left( 0, 2+\e \right]$,  
satisfy conditions~(1) through (6).    
\end{proof}

Some notation will be helpful in the next  construction.

\begin{remark}
\label{r:hlb} 
Let $\X$ be a Banach space containing  an isomorphic copy of $\ell_1$.

Recall \cite{Pel,H2} that $\X$ contains an isomorphic copy of $\ell_1$ 
if and only if  ~ $\X^*$ contains   an isomorphic copy of ~$L_1$. 
Thus    $\X^*$ also contains an isomorphic copy of $\ell_2$.
An isomorphic copy of $L_1$ (resp.~$\ell_2$) in $\X^*$ 
will be denoted by $\Z_1$ (resp.~$\Z_2$).  

There is a    norm 
$\elnm \cdot \ernm$ on $\Z_2$  
which is equivalent to the usual norm on~$\X^*$ 
and for which   $\left( \Z_2, \elnm \cdot \ernm \right)$ 
is Hilbertian; 
$\wt  \Z_2$ denotes 
$  \Z_2$ equipped with the 
new $\elnm \cdot \ernm$-norm.
Since $\wt \Z_2$ is isometric to a Hilbert space, 
there is a unique inner product that induces its 
$\elnm \cdot \ernm$-norm; 
in $\Z_2$,   Hilbert space concepts are understood 
to be in $\wt \Z_2$.     
For example, a subset of $ \Z_2 $  
is  {\it orthonormal\/}     
if, when viewed as a subset of $\wt \Z_2$, 
 it is orthonormal in $\wt \Z_2$.  
A sequence   $\{ \Y_i \}$ of finite-dimensional subspaces of $\Z_2$ 
is an {\it orthogonal finite-dimensional 
decomposition\/} ($\bot$-fdd) provided  
$\Y_i \bot \Y_j$ for $i \neq j$ and 
each $\Y_i$  is finite dimensional.
$\oc{\Y}{\Z_2}$ denotes the 
orthogonal complement 
of a subspace $\Y$  in~$\Z_2$.   
\qed  
\end{remark}

\begin{lemma} 
\label{l:tzo} 
Let $\X$ be a separable Banach space containing 
an isomorphic copy of $\ell_1$ and $\e >0$.
Then  $\Z_1$   can be taken so that a countable 
subset of it   $\nc$-norms~$\X$. 
\end{lemma} 

\begin{proof}
By \cite{H, DRT}  there is a $(1+\e)$-isomorphic copy of 
$L_1$ in $\X^*$ and so there is  an embedding 
$T \colon  \ell_1 \oplus_1 L_1 \hookrightarrow \X^*$  
satisfying,  for each $z \in \ell_1 \oplus_1 L_1$, 
$$
\lnm  ~  z ~\rnm_{\ell_1 \oplus_1 L_1} 
~\leq~ 
\lnm  ~ T z~\rnm_{\X^*} 
~\leq~ 
(1 +\e) \, \lnm ~ z ~\rnm_{\ell_1 \oplus_1 L_1}  \ . 
$$
Moreover, the image $\{ T \delta_n \}$ of the unit  vector basis 
of~$\ell_1$ can be assumed  to be weak$^*$-null 
(since $\X$ is separable, $\{ T \delta_n \}$ has a 
weak$^*$-convergent subsequence $\{ T \delta_{k_n} \}$, 
so just replace  $ \delta_n$  by 
$\frac12 (  \delta_{k_{2n}} -  \delta_{k_{2n+1}})$\,). 
Find a sequence $\{ x^*_n \}_{n=1}^\infty$ in $S(\X^*)$ such that 
$\{ x^*_n \}_{n=N}^\infty$ norms $\X $ for each ~$N \in \mathbb N$. 

Fix $\beta \in ( 0 , 1)$ and let 
\begin{gather*} 
   y_n^* = T \delta_n + \beta x^*_n   
\intertext{and} 
\Z_1 \deq \left[~ \left\{ y_n^* \hsp{3}\colon\hsp{3} n\in\mathbb N \right\}
 ~\cup~ TL_1~ \right] \ .
\end{gather*} 
Note that for each $n \in \mathbb N$, 
\begin{equation}
\label{eq:tzo1}
  \lnm y_n^* \rnm \hsp{3}\leq\hsp{3} 
1 + \e + \beta \ . 
\end{equation}

The operator $S \colon \ell_1 \oplus_1 L_1 \to \Z_1$ defined by 
$$
S\left( \sum \alpha_n \delta_n  \oplus_1  f \right) ~=~ 
        \sum \alpha_n  y_n^* ~+~ Tf 
$$ 
illustrates   that $\Z_1$ is 
isomorphic to 
$\ell_1 \oplus_1 L_1$.  
Indeed, fix  
$$\sum \alpha_n \delta_n  \oplus_1  f 
\hsp{5}\in\hsp{5} \spn{\delta_n}_{n\in\mathbb N} \oplus_1 L_1 \ . 
$$  
Then 
\begin{align*}
\lnm   S\left( \sum \alpha_n \delta_n \oplus_1 f \right)  \rnm_{\X^*} 
\hsp{3}&\leq\hsp{3} 
  \left( 1+\e+\beta\right)\, \sum \lav \alpha_n \rav   
~+~ \left( 1+\e \right) \, \lnm f \rnm_{L_1} \\
&\leq\hsp{3}
\left(1+\e + \beta \right) \, \lnm 
\sum \alpha_n \delta_n  \oplus_1  f \rnm_{\ell_1 \oplus_1 L_1} \ . 
\end{align*} 
On the other hand, 
\begin{align*}
\lnm   S\left( \sum \alpha_n \delta_n  \oplus_1  f \right)  \rnm_{\X^*} 
\hsp{3}&=\hsp{3} 
\lnm  T \left( \sum \alpha_n\delta_n \oplus_1 f \right) 
+ \beta\sum\alpha_n x_n^*\rnm_{\X^*} \\
&\geq\hsp{3}
\lnm \sum \alpha_n\delta_n \oplus_1  f \rnm_{\ell_1 \oplus_1 L_1} ~ - ~ 
\beta \lnm \sum\alpha_n  \delta_n\rnm_{\ell_1}  \\ 
&\geq\hsp{3}
\left( 1 - \beta \right) \,  \lnm 
\sum \alpha_n \delta_n  \oplus_1  f \rnm_{\ell_1 \oplus_1 L_1} \ . 
\end{align*}  
Thus  $\Z_1$ is 
$\left(\frac{1+ \beta + \e}{1-\beta}\right)\,$-isomorphic to 
$\ell_1 \oplus_1 L_1$, which is $\lc$-isomorphic to ~$L_1$. 

To see that $\{ y_n^* \}_{n\in\mathbb N}$ is
$\left( 1 + \frac{1+\e}{\beta} \right)$-norming for $\X$, 
fix  $x\in S(\X)$. 
Let  $\delta >0$ be such that $\delta(1+\delta) < \beta$  
and find    $n \in \mathbb N$ such that  
\begin{equation}
\label{eq:tzo2}
\lav (T\delta_n) (x) \rav \leq \delta 
\text{\hskip 30 pt and \hskip 30 pt} 
1 \leq (1+\delta) x^*_n (x) \ , 
\end{equation}
for then by \eqref{eq:tzo1} and \eqref{eq:tzo2}  
\begin{equation*} 
\frac{y_n^* (x)}{\lnm y_n^*\rnm}  
\hsp{3}\geq\hsp{3}
\frac{1}{1+\e+\beta}  
~\left(\frac{\beta}{1+\delta} ~-~ \delta\right) ~ . 
\end{equation*} 
Thus 
\[ 
   \sup_{n\in \mathbb N} ~ \frac{y_n^*(x)}{\lnm y_n^*\rnm}   
\hsp{3}\geq\hsp{3} 
\frac{\beta}{\beta + 1 + \e} \ . 
\]   
So,  for $\beta$   sufficiently close to one, 
$\{ y_n^* \}_{n\in\mathbb N}$  is $(2 + 2\e)$-norming for $\X$.  
\end{proof} 

\hsp{0} From Lemma \ref{l:tzo} and Theorem~\ref{t:celle} 
it easily follows, 
when~$\X$ is separable, 
that~(a) implies~(c) in~Theorem~\ref{t:elle}. 

\begin{theorem}
\label{t:celle} 
Let $\X$ be a separable Banach space 
containing $\ell_1$. From  Remark~\ref{r:hlb},  
let $\Z_1$ be total and $\Z_2 \subset \Z_1$. 
Let    $\bio{ a_n}{ b_n^*}_{n=1}^\infty$  
be in $\BIO{\X}{\Z_1}$ and       
fix ~ $\e, \eta >0$.       
Then there exists a  
 $[  (1+\sqrt{2}) + \e]$-bounded $wc_0^*$-stable biorthogonal system 
$\bio{x_n}{x_n^*}_{n=1}^\infty$ in $\BIO{\X}{\X^*}$   
so that 
\begin{enumerate} 
\item[\rm{(\ref{t:celle}a)}]
 $\lsp a_n \rsp_{n=1}^\infty \subset \lsp x_n \rsp_{n=1}^\infty  $ 
\item[\rm{(\ref{t:celle}b)}]  $\lsp b_n^* \rsp_{n=1}^\infty 
       \subset \lsp x_n^* \rsp_{n=1}^\infty \subset \Z_1$  
\item[\rm{(\ref{t:celle}c)}]   
$\sup_{n\in\mathbb N} \, d(x^*_n, \Z_2) < \eta$ . 
\end{enumerate} 
\end{theorem}
 
\noindent
In Section \ref{S:lowerbnd} it is shown that 
 the $[  (1+\sqrt{2}) + \e]$ 
can not be replaced with $(1+\e)$ in Theorem~\ref{t:celle}.   
The following fact helps  with the bound of the 
system in Theorem~\ref{t:celle}. 
It is due to Dvoretzky~\cite{Dv} 
and Milman~\cite{Mil}; a proof may be found in ~\cite{P}. 

\begin{fact}
\label{f:ldm} 
Let $n, m, N$ be positive integers and $\delta >0$.  
Then there is a positive integer $K = K( n, m, N,\delta)$ so that 
if 
\begin{enumerate} 
\item[\rm{(1)}] 
$Y$ is a Banach space with $K \leq \text{dim\,}Y \leq \infty$ 
\item[\rm{(2)}] 
$E$ is a $n$-dimensional subspace of $Y$ 
\item[\rm{(3)}] 
$H$ is a $m$-codimensional subspace of $Y$ 
\end{enumerate} 
then there is a subspace  $F$ of $H$ 
which is $\left( 1+\delta \right)$-isomorphic to $\ell_2^N$  
and a projection $P$ from  $E + F$ onto $F$ 
with $\textrm{ker\,} P = E$ and $\lnm P \rnm < 1+\delta$.  
\end{fact} \noindent 
In fact, they  showed that 
$P$ can be taken so that $\lnm P - I\vert_{E+ F} \rnm <  1+\delta$.

\begin{proof}[Proof of Theorem~\ref{t:celle}] 
The proof of Theorem~\ref{t:celle} is  similar to the 
proof of Theorem~\ref{t:csch}; thus,  
notation from the proof of  Theorem~\ref{t:csch} 
will be retained. 

Fix a strictly decreasing sequence $\{ \eta_k \}_{k=1}^\infty$   
converging to zero with $\eta_1 < \eta$.   
It suffices to  construct a system 
$\bio{x_n}{x_n^*}_{n=1}^\infty$ in $\BIO{\X}{\X^*}$  
along with 
a blocking~$\{ J_k \}_{k=1}^\infty$ of $\mathbb N$ 
and a sequence $\{ u_n^* \}_{n=1}^\infty$ from $\Z_2$  
satisfying   
\begin{enumerate}
\item[\rm{(1)}]
 $x_m^*(x_n) = \delta_{mn}$ 
\item[\rm{(2)}]
 $\lnm x_n \rnm ~\leq~ \left(1+\sqrt{2} \right) + \e$
\item[\rm{(3)}]
 $\lnm x_n^* \rnm ~\leq~  1  + \e $
\item[\rm{(4)}]
 $ \lnm x^*_n - u^*_n \rnm \leq \eta_k $ if $n \in J_k$ 
\item[\rm{(5)}]
 $\lsp u^*_n \rsp_{n\in J_{k_1}}$ is orthogonal to 
       $\lsp u^*_n \rsp_{n\in J_{k_2}}$  for $k_1 \neq k_2$   
\item[\rm{(6)}]
 $\lsp a_n \rsp_{n=1}^{\infty} 
       \subset \lsp x_n \rsp_{n=1}^{\infty}$
\item[\rm{(7)}]
  $\lsp b^*_n \rsp_{n=1}^{\infty} 
       \subset \lsp x_n^* \rsp_{n=1}^{\infty} \subset \Z_1 $   \ . 
\end{enumerate}   
Note that conditions (3) through (5) imply that 
 $\{x^*_n\}_{n\in \mathbb N}$  is weakly-null in~$\X^*$.
Clearly all that remains at this point is to 
show that the  $wc_0^*$-biorthogonal system 
$\bio{x_n}{x_n^*}_{n=1}^\infty$ is indeed stable, which is 
done in the  last step by using the  condition  
that $\lsp x_n^* \rsp_{n=1}^{\infty}$ stays inside of $\Z_1 $. 

Let 
$$
\delta ~\deq~ \frac{\e}{ 2 ~+~\sqrt{2} } \ . 
$$
The construction will inductively produce 
blocks   
\[
\bio{x_n}{x_n^*}_{n\in J_k}  
\textrm{\hsp{30}and\hsp{30}}  \{ u_n^* \}_{n\in J_k}  \hsp{5} . 
\]    
Fix $k \geq 1$. 
Assume that $\{ J_j \}_{0\leq j<k}$ 
along with 
$\bio{x_n}{x_n^*}_{n  \in J_k^p }$ 
and  $\{ u_n^* \}_{n  \in J^p_k }$  
have been constructed to satisfy conditions~(1) through (5).   
Now to construct $J_k$ 
along with $\bio{x_n}{x_n^*}_{n \in J_k }$ 
and   $\{ u_n^* \}_{n  \in J_k }$. 

The idea is to    find a biorthogonal system 
$\bio{z_n}{z_n^*}_{n\in J_k}$ in $\BIO{\mathcal  P_k}{\mathcal  Q_k}$   
by first finding   $\{z_{1+n_k}, z^*_{1+n_k} \}$  that helps 
guarantee condition~(6) if $k$ is odd  and     
condition~(7) if $k$   even; 
however,   $\{z_{1+n_k}, z^*_{1+n_k} \}$    
would  not necessarily satisfy conditions~(2) and ~(3)  
and   $z^*_{1+n_k}$  may be far from $\Z_2$  
and so $J_k^o$ and 
$$
    \bio{z_n}{z_n^*}_{n \in J_k^o}   
$$ 
are then constructed and   the  Haar matrix is applied to 
$\bio{z_n}{z_n^*}_{n\in J_k}$ to produce $\bio{x_n}{x_n^*}_{n\in J_k}$  
 and $\{ u_n^* \}_{n \in J_k }$ so that   
$$
  \bio{x_n}{x_n^*}_{n\in  J^p_k \, \cup J_k }
  \text{\hskip 25 pt and \hskip 25 pt}
  \{ u_n^* \}_{n \in  J^p_k \, \cup J_k } 
$$
satisfy conditions~(1) through (5).

Find $\{z_{1+n_k}, z^*_{1+n_k} \}$ 
just as in the proof of Theorem~\ref{t:csch}:  
in the case that ~$k$ is odd, 
be sure to choose $y^*_{1+n_k}$ in $\Z_1$, 
which is possible since $\Z_1$ is total.

Find a natural number $m_k$ larger than one   so that 
$$
2^{{-m_k}/{2}} \ 
\max\left(~ \lnm z_{1+n_k} \rnm ~,~ \lnm z_{1+n_k}^* \rnm ~\right) 
~<~ 
\min\left(~ \delta  ~,~ \eta_k  ~\right)    \ . 
$$
Let 
\begin{align*}
E_k ~&\deq~  
\left[ \{ x_n^*  \}_{n\in J_k^p} \cup  \{ z_{1+n_k}^*\} \right] \\
H_k ~&\deq~ 
   \left(~ \oc{\lsp u_n^* \rsp_{n\in J_k^p}}{\Z_2}~ \right)
    ~\cap~ 
    \left[ x_n \right]_{n\in J_k^p}^\perp  
     ~\cap~ \left[ z_{1+n_k}\right]^\perp \\
Y_k ~&\deq~ \left[ \Z_2 \cup E_k \right]  \\
N_k ~&\deq~ 2^{m_k} ~-~ 1 \hsp{3}=\hsp{3} \card{J_k^o} \ . 
\end{align*}   
Use Fact~\ref{f:ldm}  to find a subspace $F_k$ 
of $H_k$, a   projection $P_k$ with kernel $E_k$,  
and a norm one isomorphism ~$T_k$ so that 
$$
E_k + F_k 
\hsp{5} \overset{P_k}{\twoheadrightarrow}  \hsp{5}
F_k 
\hsp{5} \overset{T_k}{\longrightarrow} \hsp{5} \ell_2^{N_k} 
$$
and 
\begin{equation*}
  \max\left(~   \lnm T_k^{-1} \rnm \hsp{1},\hsp{1}  
     \lnm P_k \rnm^{2}  ~\right)   ~<~   1+\delta   \ . 
\end{equation*} 
Let $\{ e_n \}_{n\in J_k^o}$ be an orthonormal basis for $\ell_2^{N_k}$. 
For each $ n\in J_k^o$, let 
\[
z_n^* = T_k^{-1}e_n
\]
 and, using Local 
Reflexivity,  
find $z_n \in \X$ that  agrees, on $E_k$ and $F_k$,  
with a norm-preserving 
Hahn-Banach extension of $P_k^* T_k^*  e_n \in (E_k + F_k)^*$ 
to $\X^*$ and satisfies 
\[ 
\lnm z_n \rnm ~<~ \sqrt{1+\delta}\hsp{4} \lnm  P_k^* T_k^*  e_n \rnm \ . 
\]
Then   $\bio{z_n}{z_n^*}_{n\in J^o_k}$ 
is a biorthogonal system in  $\BIO{E_k^\pperp}{F_k}$ and     
$$
\max \ \left\{ ~\lnm  z_n \rnm \hsp{2}\colon\hsp{2}  {n\in J^o_k}  ~\right\}
 ~<~ 1+\delta   \ . 
$$      
Now apply the Haar matrix   
to 
$\bio{z_n}{z_n^*}_{n\in J_k}$ to 
produce  $\bio{x_n}{x_n^*}_{n\in J_k}$ 
and let 
$$
u^*_n ~\deq~  \sum_{j\in J^o_k} a^{m_k}_{nj}  z_{j}^*  
$$ 
for each $n$ in $J_k$. 

With help from the observations in Remark~\ref{r:hwm}, 
note that for each $n$ in~$J_k$ 
\begin{align*}
\lnm x_n \rnm 
~&\leq~ 
2^{{-m_k}/{2}} \   \lnm z_{1+n_k} \rnm ~+~ 
\left( 1+ \sqrt{2} \right) \max_{j\in J^o_k} \lnm z_{j} \rnm \\
~&\leq~ 
\delta ~+~ \left( 1+ \sqrt{2} \right) \left( 1+  \delta \right)\\
~&\leq~
\left(1+\sqrt{2} \right) + \e
\end{align*}
and  
\begin{align*} 
\lnm x_n^* \rnm 
~&\leq~ 
2^{{-m_k}/{2}} \   \lnm z_{1+n_k}^* \rnm ~+~ 
\lnm \sum_{j\in J^o_k} a^{m_k}_{nj}  z_{j}^* \rnm   \\
~&\leq~ 
\delta ~+~  (1+\delta)~   
\lnm \sum_{j\in J^o_k} a^{m_k}_{nj}   e_j \rnm    \\
~&\leq~
1 ~+~ 2\,\delta  \\ 
~&\leq~
1 ~+~ \e \ . 
\end{align*} 
and since \ $ x_n^* ~-~ u_n^* ~=~ a^{m_k}_{n1} \, z^*_{n_k +1}$, 
$$ 
\lnm x_n^* ~-~ u_n^* \rnm
~=~
 2^{{-m_k}/{2}}   \, \lnm z^*_{n_k +1} \rnm 
 ~\leq~ \eta_k  \ .
$$
Thus, 
$$
  \bio{x_n}{x_n^*}_{n\in  J^p_k \, \cup J_k }
  \text{\hskip 25 pt and \hskip 25 pt}
  \{ u_n^* \}_{n \in  J^p_k \, \cup J_k } 
$$
clearly satisfy conditions~(1) through (5).

This completes the inductive construction of 
the system $\bio{x_n}{x_n^*}_{n=1}^\infty$  
in $\BIO{\X}{\X^*}$, along with the blocking ~
$\{ J_k \}$ of~$\mathbb N$ and the sequence $\{ u_n^* \}_{n=1}^\infty$ 
from~ $\Z_2$, 
that satisfy   conditions~(1) through~(7).  
  
The  last step is  to verify   
that  $\bio{x_n}{x_n^*}_{n=1}^\infty$  is 
indeed stable, which, by  Lemma~\ref{l:fac}, is equivalent to 
verifying that the operator 
$S \colon \X \to c_0$ given 
by $S\left(x\right) = \left( x^*_n \left(x\right)\right)$  
factors through an injective space.     
Towards this,  
consider the following  commutative  diagram  
\vsp{5}
\[
\begin{picture}(100,50)
\put(0,0){\makebox(0,0){$\ell_1$}}
\put(100,0){\makebox(0,0){$\X^*$}}
\put(50,50){\makebox(0,0){$L_1$}}
\put(20,25){\makebox(0,0){\vector(1,1){22}}}
\put(80,25){\makebox(0,0){\vector(1,-1){22}}}
\put(50,5){\makebox(0,0){\vector(1,0){45}}} 
\put(10,30){\makebox(0,0){$A$}}
\put(86,30){\makebox(0,0){$j$}}
\put(50,10){\makebox(0,0){$B$}}
\end{picture}
\]
\vsp{5}\noindent
where $j$ is an isomorphic embedding with range $ \Z_1$  and 
$A$ and $B$ are given by
\[
A(\delta_n) = j^{-1}_o x^*_n \deq r_n 
\textrm{\hsp{15}and\hsp{15}} 
 B(\delta_n) = x^*_n \ . 
\]

Since $A^*$ and $B^*$ are of the form 
$$
A^*(f) = \left( f \left( r_n \right) \right)_n 
\text{\qquad and \qquad} 
B^*(x^{**}) = \left( x^{**}\left( x^*_n \right) \right)_n \ , 
$$ 
their  ranges are  contained in $c_0$; let 
$A^*_o$ and $B^*_o$ be the corresponding maps 
with their ranges restricted to $c_0$.  
Thus the following diagram commutes. 
\vsp{5}
\[
\begin{picture}(100,100)
\put(50,0){\makebox(0,0){$\X$}}
\put(0,50){\makebox(0,0){$\X^{**}$}}
\put(100,50){\makebox(0,0){$c_0$}}
\put(50,100){\makebox(0,0){$L_\infty$}}
\put(20,75){\makebox(0,0){\vector(1,1){22}}}
\put(80,75){\makebox(0,0){\vector(1,-1){22}}}
\put(50,55){\makebox(0,0){\vector(1,0){45}}} 
\put(10,80){\makebox(0,0){$j^*$}}
\put(86,80){\makebox(0,0){$A^*_o$}}
\put(50,60){\makebox(0,0){$B^*_o$}}
\put(20,25){\makebox(0,0){\vector(-1,1){22}}}
\put(80,25){\makebox(0,0){\vector(1,1){22}}}
\put(10,20){\makebox(0,0){$\delta$}}
\put(86,20){\makebox(0,0){$S$}}
\end{picture}
\]
\vsp{5}\noindent 
An appeal to Lemma~\ref{l:fac} finishes the proof. 
\end{proof}

Theorems~\ref{t:sch} (c) is much easier to proof if one drops the
total condition since then one can use the technique of 
Davis-Johnson-Singer  (\cite[Thm.~1]{DJ}
and~\cite[Prop.~1]{S2}).  As a partial illustration of this, we offer
the following theorem, which gives a weaker result but a smaller
constant than is provided by Theorems~\ref{t:sch} (c).

\begin{theorem}
\label{t:schf} 
Let $\X$ be a separable Banach space not containing $\ell_1$ 
such that   $\X^*$ fails the Schur property.  
Fix~$\e>0$. 
Then there is a $(2 + \e)$-bounded     
fundamental  $wc_0^*$-biorthogonal 
system $\bio{x_n}{x_n^*}$  in $\BIO{\X}{\X^*}$. 
\end{theorem}
\noindent 
The meat in the proof of Theorem~\ref{t:schf}  
is the following lemma. 

\begin{lemma}
\label{l:schf} 
Let $\X$ be a separable Banach space   
such that   $\X^*$ fails  the Schur property.   
Fix  $\e>0$.  
Then there is a $wc_0^*$-biorthogonal sequence $\bio{x_n}{ x_n^*}$   
in $\BIO{S(\X)}{\X^*}$ satisfying 
\begin{enumerate} 
\item[\rm{(1)}]  $\lnm x_n^* \rnm \leq 2+\e $   
\item[\rm{(2)}]  $\{ x_n \}$ is basic 
\item[\rm{(3)}]  $\lsp  x_n^* \rsp^\pperp + \lsp x_n \rsp$ is 
       dense in $\X$. 
\end{enumerate} 
\end{lemma}
 
\begin{proof}[Proof of Lemma~\ref{l:schf}]
Fix 
a normalized weakly-null sequence~$\{ w_n^* \}$ in~$\X^*$,  
a dense sequence $\{ d_n \}$ in~$\X$, 
a sequence $\{ \e_n \}$  decreasing to zero, 
and 
a sequence $\{ \tau_n \}$ such that $\tau_n > 1$ 
and $\Pi \tau_n < \infty$.
It is sufficient to construct 
\begin{enumerate}
\item[\rm{(a)}]  a sequence $\{ x_n \}_{n\geq 1} $ in $S(\X)$
\item[\rm{(b)}]  a sequence $\{ \wt x_n^* \}_{n\geq 1} $ in $S(\X^*)$
\item[\rm{(c)}]  finite sets $\{ F_n \}_{n\geq 0}$ in $S(\X^*)$ 
              with $F_0 \deq \emptyset$ 
\item[\rm{(d)}]  an increasing sequence $\{ k_n \}_{n=1}^\infty$ 
              of integers
\end{enumerate}
that satisfy
\begin{enumerate}
\item[\rm{(4)}] $x_n \in F_{n-1}^\pperp 
             \cap \{x^*_i\}_{i <n}^\pperp \deq \X_n$ 
\item[\rm{(5)}] $\wt x_n^* \in \{d_i\}_{i < n}^\perp 
             \cap \{x_i\}_{i <n}^\perp \deq \Y_n$ 
\item[\rm{(6)}]  $  \frac{1}{2+\e} \leq \wt x_n^* (x_n) $
\item[\rm{(7)}]  $\lnm \wt x_n^* - w_{k_n}^* \rnm \leq \e_n$
\item[\rm{(8)}] if $x \in \lsp x_j \rsp_{j\leq n}$, 
             then there is $f \in F_n$ 
              with $\lnm x \rnm \leq \tau_n f(x)$.
\end{enumerate} 
For then  just take $x^*_n = \wt x^*_n / \wt x^*_n(x_n)$.  
Note that (4) and (8) imply (2) while ~(5) and 
biorthogonality 
imply (3) since each $d_i$ has the form 
$$
d_i = \left( d_i - \sum_{n=1}^i x^*_n(d_i) x_n \right)  
  ~+~ \sum_{n=1}^i x^*_n(d_i) x_n  \ . 
$$

The construction is by induction on $n$. 
To start,  let $\wt x_1^* = w_1^*$.
Find   ~$x_1$   in ~$S(\X) $ that satisfies~(6) 
and $F_1$ that satisfies~(8). 

Fix $n > 1$ and assume that  
the items in~(a) through~(d) have been constructed up 
through the $(n-1)^{\text{th\/}}$-level.  From   
this it is possible to find $\X_n$ and~$\Y_n$. 

By Lemma~\ref{l:fcd}, there is a finite co-dimensional subspace $\Y$ of  
$\X^*$ that is $(2 + \tfrac{\e}{2})$-normed  by $\X_n$.    
Find  
$\wt x_n^* \in S \left( \Y \cap \Y_n \right)$ 
along with $k_n > k_{n-1}$ such that~(7) holds.   
Since $ \Y  $ is  $(2+ \tfrac{\e}{2})$-normed by $\X_n$, 
there is $x_n \in S(\X_n)$ such that 
$$ 1 \leq (2+\e) \wt x^*_n (x_n) \ . $$ 
Now find $F_n$ satisfying (8).    
\end{proof} 

\begin{proof}[Proof of Theorem \ref{t:schf}]   
 First find the biorthogonal 
system $\bio{x_n}{x_n^*}$ given by Lemma~\ref{l:schf}.  
 The next step is to  perturb  this system 
to produce  the desired  system.

Begin by  finding 
a bijection $p\colon \mathbb N \times \mathbb N \to \mathbb N$
satisfying  
\begin{enumerate}
\item[\rm{(i)}]   $\{ p(n,i)\}_{i=1}^\infty$ is an 
             increasing sequence 
\end{enumerate}
for each $n$ in $   \mathbb N$.
Take a dense set $\{ y_n \}$ 
in $B\left( \lsp x^*_n \rsp^\pperp\right)$. 
The underlying idea is to  use $\{ x_{p(n,i)} \}_i$ 
to capture $y_n$,  along with  $\{ x_{p(n,i)} \}_i$,
in the   span of a   small perturbation of $\{ x_{p(n,i)} \}_i$. 
 
Towards this, with the help of  (2) 
and the fact that  $\{ x_{p(n,i)} \}_i$  is not equivalent 
              to the unit vector basis of $\ell_1$, 
for each $n$ find  a sequence $\{ a_{p(n,i)} \}_i$ 
such that 
\begin{enumerate}
\item[\rm{(ii)}] $\sum_i \lav a_{p(n,i)} \rav = \infty$
\item[\rm{(iii)}] $ \sum_i a_{p(n,i)} x_{p(n,i)} \in \X$. 
\end{enumerate} 
Let 
$$
  w_{p(n,i)} \deq x_{p(n,i)} ~-~  
              \e \, \left( \sign a_{ p(n,i)} \right) \,   y_n    \ . 
$$  
Clearly $\bio{w_n}{x_n^*}$ is a 
$\left[(1 + \e)  (2 + \e) \right]$-bounded 
$wc_0^*$-biorthogonal system. 

Fix  $n_0\in \mathbb N$ 
and consider $x^* \in \lsp w_{ p(n_0,i)} \rsp_i^\perp$. 
For each   $m\in \mathbb N$,
$$
x^* \left(\sum_{i=1}^m  
                a_{ p(n_0,i)} x_{ p(n_0,i)} \right) 
 = \e x^* (y_{n_0}) 
      \sum_{i=1}^m 
       \lav a_{p(n_0,i)} \rav \ . 
$$
Combined with (ii) and  (iii),  
this gives that $x^* \in \lsp y_{n_0} \rsp^\perp$,  
which in turn implies that  
$x^* \in \lsp x_{ p(n_0,i)} \rsp_i^\perp$. 
Thus 
$$ 
\lsp ~\left\{ x_{ p(n_0,i)}  \right\}_i  \cup \left\{ y_{n_0} \right\} ~\rsp 
~\subset ~\lsp w_{ p(n_0,i)} \rsp_i \ . 
$$ 
Combined with  ~(3), 
 it follows that  $\{ w_n \}_{n\in\mathbb N}$  is fundamental. 
\end{proof}


\section{BOUNDED FUNDAMENTAL BIORTHOGONAL SYSTEMS}
\label{S:lowerbnd}

The knowledgeable reader notices that, 
in our proof of Theorem~\ref{t:celle}, 
a combination of the~\cite{OP}--method 
(which produces   $[(1+\sqrt{2})^2 + \e]$-bounded systems) 
and the~\cite{P}--method 
(which produces  $(1+\e)$-bounded systems)  is used 
to produce a   $(1+\sqrt{2}  + \e)$-bounded system.  
Using just the ~\cite{P}--method in our proof of Theorem~\ref{t:celle}  
will not guarantee that  the  
 $ x^*_n $'s are in  $\Z_1$ nor close to $\Z_2$.  
This difficulty is  not purely technical. 
For indeed, consider the below special case of Theorem~\ref{t:celle}.   

\begin{corollary}  
\label{cit}
Fix $\e, \eta > 0$.  
Let $\Z_1$ be a total subspace of $\ell_\infty=\ell_1^*$
that is isomorphic to $L_1$ 
and $\Z_2$ be a subspace of $\Z_1$ that is $(1+\eta)$-isomorphic to $\ell_2$. 
Then there exists a  
 $[  (1+\sqrt{2}) + \e]$-bounded fundamental biorthogonal system 
$\bio{x_n}{x_n^*}_{n=1}^\infty$ in $\BIO{S(\ell_1)}{\ell_\infty}$ 
satisfying 
\[ 
\sup_{n\in\mathbb N} \ d\left(x^*_n , \Z_2 \right) < \eta \ . 
\]
\end{corollary} \noindent 
Lemma~\ref{l:tz} shows that such subspaces $\Z_1$ and $\Z_2$ 
in Corollary~\ref{cit} do exist. 
Corollary~\ref{it} shows that in Corollary~\ref{cit}, the 
$ [  (1+\sqrt{2}) + \e ]$ can not be replaced with~$(1+\e)$.  
However, if the requirement~(\ref{t:celle}a) in the statement 
of Theorem~\ref{t:celle} is removed (which would basically 
remove the fundamental condition), then 
the \cite{P}--method can be used to obtain 
this variant of Theorem~\ref{t:celle} 
with $(1+\e)$ replacing $ [  (1+\sqrt{2}) + \e ]$. 

The proof of Lemma~\ref{l:tzo} gives that  for each $\e>0$ there 
exists 
\begin{align*} 
\textrm{a strictly increasing surjective function\hsp{10}} 
i &\colon\hsp{3} (0,1) \to (3(1+\e) , \infty) \\
\textrm{a strictly decreasing surjective function\hsp{10}} 
n &\colon\hsp{3} (0,1) \to (2+ \e , \infty)
\end{align*} 
so that if $\X$ is a separable Banach space 
whose dual contains an isomorphic copy of 
 $L_1$, then  for each $\beta \in (0,1)$,  
there is a subspace $\Z_1$  of $\X^*$ 
which   is ~$i(\beta)$--isomorphic to $L_1$ 
and which has a countable subset that $n(\beta)$--norms ~$\X$. 
However, if the dual space contains an isometric copy $\Z$ of $L_1$  
then the above isomorphism constant~ $i(\beta)$ can be improved.

\begin{lemma} 
\label{l:tz} 
There exists 
\begin{align*} 
\textrm{a strictly increasing surjective function\hsp{10}} 
i &\colon\hsp{3} (0,1) \to (1 , \infty) \\
\textrm{a strictly decreasing surjective function\hsp{10}} 
n &\colon\hsp{3} (0,1) \to (2  , \infty)
\end{align*} 
so that if $\X$ is a separable Banach space 
whose dual contains an isometric copy of $\ell_1$ that is 
contractively complemented in some subspace   $\Z$ of $\X^*$, 
then  for each $\beta \in (0,1)$,  
there is a subspace $\Z_1$  of $\X^*$ 
which   is $i(\beta)$--isomorphic to $\Z$
and which has a countable subset that $n(\beta)$--norms $\X$.  
\end{lemma} 
\noindent 
Since $\ell_\infty$ contains an 
isometric copy $\Z$ of $L_1$, which in turn 
contains a contractively 
complemented subspace which is isometric to $\ell_1$, 
for each positive~$\eta $, 
applying Lemma~\ref{l:tz} with $\beta$ sufficiently close to zero 
gives that there is a total subspace $\Z_1$ 
 of $\ell_\infty$ that is $(1+\eta)$--isomorphic to $L_1$, 
which in turn contains a subspace $\Z_2$ which is   $(1+\eta)$--isomorphic 
to $\ell_2$.

\begin{proof}[Proof of Lemma \ref{l:tz}] 
Find $\{ e_n^*\}_{n\in\mathbb N}$ in $\Z$ which is 
$1$--equivalent to the standard unit vector basis of~$\ell_1$ 
and a surjective contractive projection 
\[
 P \colon \Z \to[e^*_n]_{n\in \mathbb N} \ . 
\] 
Without loss of generality, $\{ e_n^*\}_{n\in\mathbb N}$ 
is   weak$^*$-null 
(similar to before, just replace $\{e^*_n\}$ with a 
weak$^*$-convergent 
subsequence $\{\frac{1}{2}(e^*_{k_{2n}} - e^*_{k_{2n+1}})\}$, 
which will be $1$--equivalent to $\{ e^*_n\}$ and contractively complemented
in $[ e^*_n ]$).    
Find a sequence~$\{ x^*_n \}_{n=1}^\infty$ in $S(\X^*)$ such that 
$\{ x^*_n \}_{n=N}^\infty$ norms $\X $ for each ~$N \in \mathbb N$. 

Fix $\beta \in ( 0 , 1)$ and let 
\begin{gather*} 
   y_n^* = e_n^* + \beta x^*_n   
\intertext{and} 
\Z_1 \deq \left[~ \left\{ y_n^* \hsp{3}\colon\hsp{3} n\in\mathbb N \right\}
 ~\cup~  \textrm{ker } P~ \right] \ .
\end{gather*} 
Note that for each $n \in \mathbb N$, 
\begin{equation}
\label{eq:tz1}
  \lnm y_n^* \rnm \hsp{3}\leq\hsp{3} 
1 +   \beta \ . 
\end{equation} 

Each element in $\Z$ has a unique expression as 
\[ 
  z = \sum \alpha_n e_n^* + f 
\] 
where $f\in \textrm{ker }P$; the 
operator $S \colon \Z \to \Z_1$ defined by 
$$
S\left( \sum \alpha_n  e_n^* +  f \right) ~=~ 
        \sum \alpha_n  y_n^* ~+~  f 
$$ 
illustrates   that $\Z_1$ is 
 isomorphic to 
$\Z$.  
Indeed, fix  
$$\sum \alpha_n e_n^*  +  f 
\hsp{5}\in\hsp{5} \spn{e_n^*}_{n\in\mathbb N} + \textrm{ker }P  \ . 
$$  
Then 
\begin{align*}
\lnm z - Sz \rnm 
\hsp{3}&=\hsp{3} 
\lnm \sum \alpha_n (e^*_n - y^*_n ) \rnm \\
&\leq\hsp{3}
\beta \sum \lav \alpha_n \rav  
\hsp{5}= \hsp{5}  
\beta \lnm P(z) \rnm  \\ 
&\leq\hsp{3}
\beta \lnm z \rnm \ . 
\end{align*}
Thus 
\[
(1-\beta) \, \lnm z \rnm   \hsp{5}\leq\hsp{5} 
\lnm Sz \rnm   \hsp{5}\leq\hsp{5}
(1+\beta) \, \lnm z \rnm 
\]
for each $z\in \Z$; thus,   
$\Z_1$ is $(\frac{1+ \beta}{1-\beta})\,$-isomorphic to $\Z$. 

To see that $\{ y_n^* \}_{n\in\mathbb N}$ is
$\left( 1 + \frac{1}{\beta} \right)$-norming for $\X$, 
fix  $x\in S(\X)$. 
Let  $\delta >0$ be such that $\delta(1+\delta) < \beta$  
and find    $n \in \mathbb N$ such that  
\begin{equation}
\label{eq:tz2}
\lav  e^*_n  (x) \rav \leq \delta 
\text{\hskip 30 pt and \hskip 30 pt} 
1 \leq (1+\delta) x^*_n (x) \ , 
\end{equation}
for then by \eqref{eq:tz1} and \eqref{eq:tz2}  
\begin{equation*} 
\frac{y_n^* (x)}{\lnm y_n^*\rnm}  
\hsp{3}\geq\hsp{3}
\frac{1}{1+ \beta}  
~\left(\frac{\beta}{1+\delta} ~-~ \delta\right) ~ . 
\end{equation*} 
Thus 
\[ 
   \sup_{n\in \mathbb N} ~ \frac{y_n^*(x)}{\lnm y_n^*\rnm}   
\hsp{3}\geq\hsp{3} 
\frac{\beta}{\beta + 1  } \ . 
\]   
So $\{ y_n^* \}_{n\in\mathbb N}$ 
is $(1 + \frac{1}{\beta})$-norming for $\X$.  
\end{proof}  
  
Recall that the modulus of  convexity 
$\mocw  \colon [0,2] \to [0,1]$ of a Banach space~$\W$ is 
$$
\mocw (\e) \deq 
\inf \left\{ 1 - \lnm \frac{x+y}{2} \rnm ~\colon~ \hsp{3}
x,y \in B(\W)  \text{\hskip 5 pt and \hskip 5 pt} 
\lnm x - y \rnm \geq \e \right\}         
$$  
and  $\W$ is {\it uniformly convex\/} if    
$\moc_\W (\e)>0$ for each $\e\in (0,2]$.    
If ~$\W$ is uniformly convex   
then $\moc_\W$ is a surjective continuous strictly increasing   
function (cf. ~\cite[pp.~53--55]{GK}).

In a uniformly convex space, 
the midpoint of points  near to  the sphere that are far apart 
is uniformly bounded away from the sphere.  
This can be extended to convex combinations of points 
near to the sphere.

\begin{lemma}
\label{l:uc} 
Let $\W$ be a   uniformly convex Banach space 
and 
$\e $ and $b$ be constants satisfying   
$$
0 ~<~ \e ~\leq b ~\leq~  1 \ . 
$$ 
If $\{ \alpha_j \}_{j=1}^\infty \subset  \mathbb  R$ 
and $\{ y_j \}_{j=1}^\infty \subset B(\W)$
satisfy       
\begin{equation}
\label{eq:lucc}
\sum_{j=1}^\infty \lav \alpha_j\rav = 1
\text{\hsp{20}and\hsp{25}} 
\lnm  \sum_{j=1}^\infty \alpha_j y_j   \rnm ~>~ 1 ~-~   \e \ , 
\end{equation}
then there is a finite subset $F$  of $\mathbb  N$ so that 
\begin{equation}
\label{eq:luclb}
\sum_{j\in F} \lav \alpha_j \rav ~>~ 1~-~ \frac{\e}{2b - \e} 
\end{equation} 
and for each $j  \in F$, 
\begin{equation}
\label{eq:lucub}
\lnm (\sign \alpha_j)\, y_j ~-~   \sum_{j=1}^\infty \alpha_j y_j 
   \rnm ~<~ \moc^{-1}_\W \left( b \right)  
\end{equation} 
and $ \alpha_j \not = 0 $  . 
\end{lemma}

\begin{proof} 
Let $\{ \alpha_j \}_{j=1}^\infty \subset  \mathbb  R$ 
and $\{ y_j \}_{j=1}^\infty \subset B(\W)$  satisfy~\eqref{eq:lucc}. 
Set 
$$
x_0 \deq \sum_{j=1}^\infty \alpha_j y_j \ . 
$$ 
Without loss of generality  each $\alpha_j  \geq 0$. 

Find $x_0^* \in S(\W^*)$ so that $x^*_0(x_0) = \lnm x_0 \rnm$ 
and let 
$$ 
F = \{ j \in \mathbb  N ~\colon~ x^*_0(y_j) ~>~ 1 ~+~ \e ~-~ 2b 
 \text{\hsp{7} and \hsp{7}} \alpha_j \not = 0 \}. 
$$ 
The condition $\e \leq b$ guarantees that 
$F$ is non-empty. 
Since  
\begin{align*}  
1-    \e 
~&<~ 
\sum_{j\not\in F} \alpha_j  x^*_0(y_j) ~+~ \sum_{j \in F} \alpha_j  x^*_0(y_j) \\
~&\leq~   
\left(1 + \e - 2b\right)
   \left(1 ~-~ \sum_{j\in F} \alpha_j\right)
~+~ \left(\sum_{j \in F} \alpha_j\right) \\
&=~ 
\left( 1 + \e - 2b  \right)
+  \left( 2b - \e \right)~\left(\sum_{j \in F} \alpha_j\right)  \ , 
\end{align*}
condition \eqref{eq:luclb} holds. 
For  each $j\in F$,
$$
\frac12 \lnm  y_j + x_0 \rnm 
~>~ 
\frac12 
\left( ~  \left( 1 + \e -2b \right) ~+~ 
   \left(1 - \e \right) ~ \right)  
~=~ 
1 - b    
$$ 
and so     
$$
\lnm y_j  ~-~  x_0  \rnm ~<~ 
\moc^{-1}_\W\left( b \right)  
$$
by uniform convexity. 
\end{proof}

\begin{proposition}
\label{lbeq}
Let $\bio{x_n}{x_n^*}_{n=1}^\infty$ be a    biorthogonal system 
in $\BIO{S\left(\ell_1\right)}{\ell_\infty}$     
and~$\W$   be a  uniformly convex Banach space 
and  
$Q~\in~\mathcal  L\left(\ell_1, \W\right) $  be  
of norm at most one, all of which satisfy  
\begin{equation}
\label{eq:lbce}
d\left( x_n^* , Q^* \left( K\, B\left(\W^*\right)\right) \right) ~<~  2 \eta 
\end{equation} 
for some constants $K\ge 1$ and $\eta > 0$. 
If 
\begin{equation}
\label{eq:lbc}
1 ~-~ \frac{1-2\eta}{K} 
~<~
\frac{2}{3}~  
\moc_\W\left( \frac{1-4\eta}{2K}\right) 
\end{equation} 
then $\{x_n\}_{n=1}^\infty$ is equivalent
to the standard unit vector basis  of~$\ell_1$. 
More specifically 
if constants $a$ and $ b$  
satisfy  
\begin{equation} 
\label{eq:lbab}
\frac{3}{2}\left(1 ~-~ \frac{1-2\eta}{K}\right)    
~\leq~ a   ~<~ b 
~\leq 
\moc_\W\left( \frac{1-4\eta}{2K}\right)  
\end{equation}   
then  
\begin{equation} 
\label{eq:lbeq}
\frac{ 3 \left(b-a\right)}{3b-a} ~
\sum_{n=1}^\infty \lav \beta_n\rav
~\leq~ 
\lnm \sum_{n=1}^\infty   \beta_n x_n \rnm 
~\leq \sum_{n=1}^\infty \lav \beta_n\rav 
\end{equation}
for each $\{ \beta_n \}_{n=1}^\infty$ in $\ell_1$.  
\end{proposition}
  
\noindent
Note 
that if $K = 1$ and $\eta = 0$ then \eqref{eq:lbc} becomes 
$$ 
 0 ~<~ \frac{2}{3}~ 
\moc_\W \left( \frac{1}{2}  \right) \ ;  
$$ 
thus, if $K$ is sufficiently close to $1$ 
and    $\eta$ is sufficiently close to   $0$, 
as they often are in  practice, then   
\eqref{eq:lbc} does indeed hold.

\begin{proof} 
The underlying idea behind the proof is to use  Lemma~\ref{l:uc} 
to find a small per\-tur\-ba\-tion 
$\{ \wt{x}_n \}_{n=1}^\infty$ 
of    $\{x_n\}_{n=1}^\infty$  that are disjointly supported  
on the standard unit vector basis  of~$\ell_1$. 
For then,  
$\{ \wt{x}_n \}_{n=1}^\infty$  is equivalent 
to the 
standard unit vector basis of $\ell_1$ and so, for a 
small enough perturbation,~$\{  x_n \}_{n=1}^\infty  $  
is also equivalent to the 
standard unit vector basis of $\ell_1$. 

Find $\{ w^*_n \}_{n=1}^\infty$ in $K\, B\left(   \W^*\right)$   so that  
\begin{equation}
\label{eq:dist}
\lnm   x^*_n  ~-~ Q^* w^*_n \rnm ~<~  2\eta    \ . 
\end{equation}   
Thus  
\begin{align*} 
\lnm Q x_n \rnm  
&\geq~
\lav \left\langle Q x_n ,  
  \frac{w^*_n}{\lnm w^*_n\rnm} \right\rangle \rav \\ 
&\geq~
\frac{1}{\lnm   w^*_n\rnm}  ~
\lav \left\langle x_n , x^*_n \right\rangle 
~-~ 
\left\langle x_n ,  x_n^* - Q^* w^*_n  \right\rangle \rav \\ 
&>~
\frac{1-  2\eta}{K}   
\end{align*}
and so by   the   first inequality in \eqref{eq:lbab}   
$$ 
\lnm Q x_n \rnm ~>~ 1 ~-~  2a/3 \ . 
$$    
 
Write 
$$
x_n = \sum_{j=1}^\infty ~ \alpha_j^n \delta_j 
$$
where 
$\sum_{j=1}^\infty ~ \lav \alpha_j^n\rav = 1$ and 
let $\e^n_j =  \sign \alpha_j^n$; thus,   
$$ 
Q x_n = \sum_{j=1}^\infty ~ \alpha_j^n ~ Q \delta_j  \ .  
$$ 
So  by  Lemma~\ref{l:uc} 
there is a sequence  $\{F_n \}_{n=1}^\infty $ 
of finite subsets in $\mathbb  N$  so that 
$$
 \sum_{j\in F_n} \lav \alpha_j^n \rav ~>~  
1 ~-~ \frac{a}{3b-a} 
$$ 
and for each   $j  \in F_n$
$$
\lnm \e^n_j Q\delta_j ~-~   
Q x_n \rnm  ~<~  \moc_\W^{-1} \left( b
\right)  
$$  
and $\alpha^n_j \not =0 $.  Let 
$$
  \wt{x}_n ~=~ \sum_{j\in F_n}   \alpha_j^n ~\delta_j \ . 
$$

To see that the $F_n$'s are disjoint, 
suppose that there is $j_0 \in F_n \cap F_m $ for some 
distinct $n,m \in \mathbb  N$. 
Pick $\tau \in \{ -1, 1\}$ so that 
$\tau \e_{j_0}^m =  \e_{j_0}^n$. 
Then 
\begin{align*} 
\lnm Q x_n ~-~ \tau Q x_m \rnm 
~&\leq~ 
\lnm Q x_n ~-~ \e^n_{j_0} Q \delta_{j_0} \rnm ~+~ 
\lnm \tau \e^m_{j_0} Q \delta_{j_0} ~-~ \tau Q x_m \rnm \\
~&<~ 2 \moc_\W^{-1}\left( b \right) \ . 
\end{align*} 
On the other hand,   
from \eqref{eq:dist} and the third inequality of \eqref{eq:lbab}  
it follows that 
\begin{align*}
 \lnm Q x_n  ~ \pm~ Q x_m \rnm 
~ &\geq ~
\lav \left\langle Q( x_n \pm x_m ) 
   ~,~ \frac{w^*_n}{\lnm w^*_n \rnm} \right\rangle \rav \\
~&\geq~
\frac{1}{\lnm w^*_n \rnm} ~
\lav \left\langle x_n \pm x_m , x^*_n \right\rangle ~-~ 
        \left\langle x_n \pm x_m , x^*_n -  Q^* w^*_n \right\rangle \rav \\
~&> ~
\frac{1 -  4\eta}{K}  \\
&\geq~
 2 \moc_\W^{-1}\left( b \right)   \ . 
\end{align*}
A contradiction, thus the finite subsets $\{ F_n \}_{n=1}^\infty$ 
of $\mathbb  N$ are indeed disjoint.

For each $\{ \beta_n \}_{n=1}^\infty$ in $\ell_1$, 
\begin{align*} 
\sum_{n=1}^\infty \lav \beta_n\rav
~&\geq 
\lnm \sum_{n=1}^\infty \beta_n  x_n \rnm  \\ 
~&\geq ~ 
\lnm \sum_{n=1}^\infty \beta_n \wt{x}_n \rnm  ~-~
 \sum_{n=1}^\infty \lav \beta_n \rav \lnm \wt{x}_n - x_n \rnm  \\
~&\geq
\left(1- \frac{a}{3b-a} \right)~ 
\sum_{n=1}^\infty \lav \beta_n\rav ~-~ 
\left(\frac{a}{3b-a} \right)  
~ \sum_{n=1}^\infty \lav \beta_n\rav   \\
~&=~
\left(1- \frac{2a}{3b-a} \right)~ \sum_{n=1}^\infty \lav \beta_n\rav  
\end{align*}
and so \eqref{eq:lbeq} holds.  
\end{proof}

If $\W$ is uniformly convex 
then $\W^*$ is super-reflexive and 
so $\W^*$ has finite cotype; thus, 
there exists   
a cotype constant $C_q\left(\W^*\right) \geq 1 $  
for some $q \in \left[2,\infty\right)$  so that 
\begin{equation}
\label{eq:ctc}
\left( \sum_{i=1}^n \lnm w_i^* \rnm^q \right)^{1/q} 
~\leq~ C_q\left(\W^*\right) ~ \left[ ~\text{avg}\sb{\theta_i=\pm 1} 
 ~\lnm \sum_{i=1}^n \theta_i w_i^* \rnm^2 ~\right]^{1/2}
\end{equation}
for each finite sequence $\{ w_i^* \}_{i=1}^n$ in $\W^*$.

\begin{theorem}
\label{lbuc}
Let   a  uniformly convex Banach space~$\W$  
 and~$\e_0>0$ 
satisfy 
\begin{equation}
\label{eq:e0}
1 - \frac{1}{1+\e_0} ~=~
\frac{2}{3} \, \moc_\W\left(\frac{1}{2(1+\e_0)}\right)  \ . 
\end{equation} 
Let  $C_q\left(\W^*\right)$ be as in~\eqref{eq:ctc} 
and $0 < \e_1 < \e_0$.  
Then   there exists  a constant
\begin{equation}
\label{eq:eta}
\eta = \eta\left(C_q\left(\W^*\right), \moc_\W, \e_1\right) > 0
\end{equation}
so that if  
\begin{enumerate}  
\item[\rm{(\ref{lbuc}a)}]   
$\bio{x_n}{x_n^*}_{n=1}^\infty$  is a 
$(1+\e)$--bounded   fundamental  biorthogonal system 
in $\BIO{S\left(\ell_1\right)}{\ell_\infty}$
\item[(\rm{\ref{lbuc}b)}]  
$Q~\in~\mathcal  L\left(\ell_1, \W\right)$
 \item[\rm{(\ref{lbuc}c)}]
$Q^*$ is a $(1+\eta)$--isomorphic embedding 
\item[\rm{(\ref{lbuc}d)}]
$d\left( x_n^* ,  Q^* \W^* \right)   ~ <~  \eta $
\end{enumerate}
 then 
$$ 
\e > \e_1  \ . 
$$  
\end{theorem}

\noindent   
The following notation helps 
crystallize condition~\eqref{eq:e0} 
and simplify some technical arguments in the 
proof of  Theorem~\ref{lbuc}.

\noindent  
\textit{Notation.} 
Consider the functions 
$l, u \colon \left[ 0, 1/4 \right] \times \left[ 0, \infty \right) 
\times  \left[ 0, \infty \right) \to [0, 1]$ given by 
\begin{equation}
\begin{split}
\label{deflu} 
l(\eta_1, \eta_2, \e) ~&=~  
1 -  \frac{1-2\eta_1}{(1+\eta_2)(1+\e)} \\
u(\eta_1, \eta_2, \e) ~&=~  \frac{2}{3}\, 
\moc_\W\left( \frac{1-4\eta_1}{2(1+\eta_2)(1+\e)}\right) \ . 
\end{split}
\end{equation}
Note that in each variable, $l$ is a
strictly increasing continuous function 
and~$u$ is a strictly decreasing continuous function. 
Condition~\eqref{eq:e0} is equivalent to 
$$
l(0,0,\e_0) ~=~ u(0,0,\e_0)  
$$
and 
\begin{align*}
l(0,0,\cdot) \hsp{3}&\colon\hsp{3} 
 [0,\infty) 
\hsp{3}\overset{\textrm{onto}}{\longrightarrow}\hsp{3} [0,1) \\
u(0,0,\cdot) \hsp{3}&\colon\hsp{3} 
 \left[0,\infty\right) 
\hsp{3}\overset{\textrm{onto}}{\longrightarrow}\hsp{3} 
\left(0,\frac{2}{3} \, \mocw\left(\frac{1}{2}\right) \right] \ ; 
\end{align*}  
thus, for a uniformly convex space $\W$, 
 there  is a unique $\e_0 >0$ satisfying~\eqref{eq:e0}.

\begin{proof}  
The underlying idea behind the proof is that  
 for  sufficiently small~$\eta$ 
Proposition~\ref{lbeq} gives that~$\{ x_n \}$ 
is equivalent to the standard unit vector basis of $\ell_1$: 
indeed, condition~\eqref{eq:lbce} will hold and 
condition~\eqref{eq:e0}  implies~\eqref{eq:lbc}. 
Then ~$\{ x_n^* \}$  is equivalent to the standard unit vector basis of $c_0$.
But if $\e$ is small enough, then condition~(\ref{lbuc}d) cannot hold 
since $\W^*$ has finite cotype.   

Let the  hypotheses of 
Theorem~\ref{lbuc}   hold.  
 Since 
$$
l(0,0,\e_1) < l(0,0,\e_0) = u(0,0,\e_0) < u(0,0,\e_1) 
$$
there are constants $a$ and $b$ so that 
$$
l(0,0,\e_1) < a < b < u(0,0,\e_1)  \ . 
$$
Find 
$\eta = \eta\left(C_q\left(\W^*\right), \moc_\W, \e_1\right) > 0$   
sufficiently small enough so that 
$$
l(\eta, \eta, \e_1) < a < b < u(\eta, \eta, \e_1)
$$
and so that there exists $N\in\mathbb N$  satisfying  
\begin{equation}
\label{eq:yuk}
 C~ (1+\eta) 
\left[ 1 + \frac{2a}{3(b-a)} ~+~ 2N\eta 
 \right] ~<~ N^{1/q}  
\end{equation}    
where $C = C_q\left(\W^*\right)$. 
To see that condition~\eqref{eq:yuk} is easily 
accomplished, note that if 
$$\left[ 2C\left(3+\frac{2a}{3(b-a)}\right)\right]^q < N \in \mathbb N$$ 
and $0 < \eta \leq \frac{1}{N}$  then~\eqref{eq:yuk} holds. 

Let  conditions (\ref{lbuc}a) through (\ref{lbuc}d) hold.  
Assume that $\e \leq \e_1$.

By (\ref{lbuc}c), 
without loss of generality, for each $w^* \in \W^*$,  
\begin{equation}
\label{eq:lbucc}
\frac{1}{1+ {\eta} }\, \lnm w^* \rnm 
~\leq~ \lnm Q^*w^* \rnm ~\leq~ \lnm w^* \rnm  \ . 
\end{equation}  
Keeping with the notation from Proposition~\ref{lbeq}, let  
$$ 
K =  (1+{\eta})\, (1+\e) \ . 
$$ From (\ref{lbuc}a), 
(\ref{lbuc}d),  and~\eqref{eq:lbucc} it follows that there 
is $\{ w_n^* \}_{n=1}^\infty \subset K B(\W^*)$ 
such that 
\begin{gather*}
\lnm x_n^* - Q^* w_n^* \rnm <  2\eta  \\
\lnm x_n^* \rnm ~=~ \lnm Q^* w^*_n \rnm \ . 
\end{gather*}
Thus~\eqref{eq:lbce} from Proposition~\ref{lbeq}   holds. 
Furthermore~\eqref{eq:lbc} also holds since  
$$
l(\eta , {\eta}, \e) \leq l(\eta , {\eta}, \e_1) 
< a < b < u(\eta , {\eta}, \e_1) \leq u(\eta , {\eta}, \e)  \ . 
$$
So by Proposition~\ref{lbeq}, since  
$\{x_n\}_{n=1}^\infty$ is   fundamental, 
$\{x_n^*\}_{n=1}^\infty$ is  equivalent
to the standard unit vector basis  of~$c_0$ with 
$$  
\frac{3\left( b-a\right)}{3b-a}
\lnm \sum_{n=1}^\infty   \beta_n x_n^* \rnm 
~\leq~ 
\sup_n \lav \beta_n \rav 
~\leq~ 
\lnm \sum_{n=1}^\infty   \beta_n x_n^* \rnm 
$$
for each $\{ \beta_n \}_{n=1}^\infty$ in $c_0$.

For each finite subset $F$ of $\mathbb  N$  
\begin{equation}
\label{eq:ct}
\left[ \sum_{n\in F} \lnm w^*_n \rnm^q 
  \right]^{\frac{1}{q}} 
~\leq~  
C  ~ \left[ ~\text{avg}\sb{\theta_i=\pm 1} 
 ~\lnm \sum_{i=1}^n \theta_i w_i^* \rnm^2 ~\right]^{1/2}
\end{equation}
The   right-hand side of~\eqref{eq:ct} 
 mimics~$c_0$-growth in that for 
each choice   $\left\{ \theta_n \right\}_{n\in F}$  of signs
\begin{align*} 
\lnm \sum_{n\in F} \theta_n   w_n^* \rnm  
~&\leq~
(1+{\eta}) ~ \lnm \sum_{n\in F} \theta_n Q^* w_n^* \rnm \\
~&\leq~
(1+{\eta})  ~\left[~ \lnm \sum_{n\in F} \theta_n  x_n^* \rnm 
+ \lnm \sum_{n\in F} 
   \theta_n \left( x^*_n - Q^* w_n^* \right) \rnm~\right] \\
~&\leq~
(1+{\eta})  ~\left[ \frac{3b-a}{3(b-a)} ~+~ 2\card{F}\eta \right]   \ .  
\end{align*}  
The    left-hand side of~\eqref{eq:ct} 
mimics~$l_q$-growth  since 
$$
1  ~\leq~ \lnm x^*_n \rnm   ~ =~ \lnm Q^*w^*_n \rnm 
~\leq~ \lnm w_n^* \rnm 
$$
for each $n\in\mathbb  N$. Thus  
$$
 \card{F}^{\frac{1}{q}} 
~\leq~ C (1+\eta)   \, 
\left[ \frac{3b-a}{3(b-a)}  ~+~ 2\card{F}\eta \right]  \ . 
$$ 
This contradicts \eqref{eq:yuk}.  Thus $\e_1 < \e$. 
\end{proof}

\begin{corollary} 
\label{it}   
Let 
$$
0 < \e_1 < - 2 + \frac{\sqrt{147}}{{6}}  
$$ 
and, following the notation in \eqref{eq:eta},  
\[
\eta = \eta\left( 1, \moc_{\ell_2} , \e_1 \right) \ . 
\] 
Let $\bio{x_n}{x_n^*}_{n=1}^\infty$  be a 
$(1+\e)$--bounded   fundamental  biorthogonal system 
in $\BIO{S\left(\ell_1\right)}{\ell_\infty}$ 
satisfying 
\[ 
\sup_{n\in\mathbb N} \ d\left(x^*_n , \Z_2 \right) < \eta 
\] 
for some subspace $\Z_2$ of $\ell_\infty$ that is a 
$(1+\eta)$--isomorph of a Hilbert space.  
Then 
\begin{equation}
\label{eq:it} 
 \e  > \e_1  \ . 
\end{equation} 
\end{corollary}

\begin{proof}  
Let $\W = \ell_2$. Thus $C_2(\W^*)=1$. 
It   is straight forward to verify that    
\[ 
\e_0 \deq - 2 + \frac{\sqrt{147}}{{6}}   
\] 
satisfies condition~\eqref{eq:e0} of Theorem~\ref{lbuc}. 
Note that $\e_0 \approx 0.0207$.

Let $\pdz$ be the predual of $\Z_2$. 
There is an operator $T~\in~\mathcal L(\ell_1, \pdz)$ 
such that $T^*\in \mathcal L (\Z_2, \ell_\infty)$ 
is the formal pointwise embedding; 
for indeed, since $\Z_2$ is reflexive, this  
formal pointwise embedding is weak$^*$-to-weak$^*$ continuous. 
Similarly, by reflexivity,   
there is $S\in \mathcal L ( \pdz, \ell_2)$ 
such that $S^*\in \mathcal L ( \ell_2 , \Z_2)$  
is a $(1+\eta)$--isomorphism. Let $Q = S \circ T$. 
Thus   
\[ 
Q \colon \ell_1 
\hsp{3}\overset{T}{\to}\hsp{3} {\pdz} 
\hsp{3}\overset{S}{\to}\hsp{3} \ell_2 
\hsp{30}\textrm{and}\hsp{30}   
Q^* \colon \ell_2 
\hsp{3}\overset{S^*}{\to}\hsp{3}  \Z_2 
\hsp{3}\overset{T^*}{\to}\hsp{3} \ell_\infty  
\]
and $Q^*$ is a $(1+\eta)$--isomorphic embedding. 

Thus condition~\eqref{eq:it} follows from
Theorem~\ref{lbuc}. 
\end{proof}


\end{document}